\documentclass[11pt]{report}
\usepackage{amsfonts}
\normalbaselines

\newcommand{\cf}[1]{\mbox{\boldmath${#1}$}}
\newcommand{\E}[0]{{\rm e}}
\newcommand{\I}[0]{{\rm i}}
\newcommand{\dis}[1]{\displaystyle{#1}}
\newcommand{\text}[1]{\textstyle{#1}}

\newcommand{\sscr}[1]{\scriptscriptstyle{#1}}

\newcommand{\vep}[0]{\varepsilon}

\usepackage{wrapfig}

\usepackage[
    colorlinks=true,
    urlcolor=black,
    filecolor=black,
    linkcolor=black,
    anchorcolor=black,
    menucolor=black,
    citecolor=black,
    pdftitle={Laguerre 2D},
    pdfauthor={Wuensche},
    pdfsubject={Laguerre 2D, article},
    pdfkeywords={Laguerre 2D, article},
    pdfproducer={pdfLaTeX},
    pdfpagelabels=true,
    hyperindex=true,
    plainpages=false,
    hypertexnames=false,
    linktocpage=true,
    pdfpagemode=UseThumbs, 
    pdfstartview=FitH, 
    pdfview=FitH,
    pdfmenubar=true,
    bookmarks=true,
    bookmarksnumbered=true,
    bookmarksopen=true,
    bookmarksopenlevel=1,
    bookmarksnumbered=true
]{hyperref}
\usepackage[pdftex]{graphicx}
\usepackage[pdftex]{thumbpdf}
\makeindex

\newcounter{saveeqn}
\newcommand{\aeqn}{\setcounter{saveeqn}{\value{equation}}%
\setcounter{equation}{0}%
\renewcommand{\theequation}{\mbox{A.\arabic{equation}}}}

\begin{document}
\bibliographystyle{unsrt}

\vbox {\vspace{0mm}}
\begin{center}
{\LARGE \bf Generating Functions for Products of Special \\[4mm] Laguerre {\bf 2D} and Hermite {\bf 2D} Polynomials}\\[3mm]
\end{center}

\begin{center}
{A. W\"{u}nsche}
\\[1mm]
formerly: {\it Institut f\"{u}r Physik, Nichtklassische Strahlung,\\ Humboldt-Universit\"{a}t,\\
Newtonstr. 15, 12489 Berlin, Germany} \\[2mm]
e-mail:$\;$ alfred.wuensche@physik.hu-berlin.de \\

\end{center}

\begin{center}

{\bf{Abstract}}

\begin{quote}

The bilinear generating function for products of two Laguerre 2D polynomials ${\rm
L}_{m,n}(z,z')$ with different arguments is calculated. It corresponds to the formula of
Mehler for the generating function of products of two Hermite polynomials. Furthermore,
the generating function for mixed products of Laguerre 2D and Hermite 2D polynomials and
for products of two Hermite 2D polynomials is calculated. A set of infinite sums over
products of two Laguerre 2D polynomials as intermediate step to the generating function
for products of Laguerre 2D polynomials is evaluated but these sums possess also proper
importance for calculations with Laguerre polynomials. With the technique of $SU(1,1)$
operator disentanglement some operator identities are derived in an appendix. They allow
to calculate convolutions of Gaussian functions combined with polynomials in one- and
two-dimensional case and are applied to evaluate the discussed generating functions.


\end{quote}

\end{center}

\setcounter{chapter}{1}
\setcounter{equation}{0}
\section*{1. Introduction}

Hermite and Laguerre polynomials play a great role in mathematics and in mathematical
physics and can be found in many monographs of Special Functions, e.g.,
[1--4]. Special comprehensive representations of polynomials of two and of several variables are given in, e.g, \cite{suetin,dunkl}.

Laguerre 2D polynomials ${\rm L}_{m,n}(z,z')$ with two, in general, independent complex
variables $z$ and $z'$ were introduced in 
[7--12] by (similar or more
general objects with other names and notations were defined in
[13--24])
\begin{eqnarray}
{\rm L}_{m,n}(z,z') &\equiv &  \exp\left(-\frac{\partial^2}{\partial z \partial
z'}\right)z^{m}z'^{\,n} \underbrace{\exp\left(\frac{\partial^2}{\partial z \partial
z'}\right)\, 1}_{=1} \nonumber\\ &=& (-1)^{m+n} \exp(zz')
\frac{\partial^{m+n}}{\partial z'^{\,m}\partial z^{n}}\exp(-z z')\, 1\nonumber\\ &=&
\left(z-\frac{\partial}{\partial z'}\right)^m \left(z'-\frac{\partial}{\partial z}\right)^n 1,
\label{lag2d}
\end{eqnarray}
written by application of an operator to the function $f(z,z')=1$. This leads to the following
definition (called "operational" in comparison to the "Rodrigues"-like) and explicit representation
\begin{eqnarray}
{\rm L}_{m,n}(z,z') &=& \exp\left(-\frac{\partial^2}{\partial z \partial
z'}\right)z^{m}z'^{\,n} \nonumber\\ &=&
\sum_{j=0}^{\{m,n\}}\frac{(-1)^j m!n!}{j!(m-j)!(n-j)!}z^{m-j}z'^{\,n-j}, \label{lag2d1}
\end{eqnarray}
with the inversion (see also formulae (\ref{difflag2d}))
\begin{eqnarray}
z^m z'^{\,n} &=& \exp\left(\frac{\partial^2}{\partial z \partial
z'}\right)\,{\rm L}_{m,n}(z,z') \nonumber\\ &=& \sum_{j=0}^{\{m,n\}}\frac{m!n!}{j!(m-j)!(n-j)!}
\,{\rm L}_{m-j,n-j}(z,z').
\end{eqnarray}
Some special cases are
\begin{eqnarray}
&& {\rm L}_{m,n}(z,0) \;=\; \frac{(-1)^n m!}{(m-n)!}z^{m-n},\quad {\rm L}_{m,n}(0,z') \;=\; \frac{(-1)^m n!}{(n-m)!}z'^{\,n-m},\quad {\rm L}_{m,n}(0,0) \;=\; (-1)^n n!\delta_{m,n}, \nonumber\\ && {\rm L}_{m,0}(z,z') \;=\; z^m,\quad {\rm L}_{0,n}(z,z') \;=\; z'^{\,n}, \quad {\rm L}_{0,0}(z,z') \;=\; 1.\label{lag2d1z0}
\end{eqnarray}
The differentiation of the Laguerre 2D polynomials provides again Laguerre 2D polynomials
\begin{eqnarray}
\frac{\partial}{\partial z}\,{\rm L}_{m,n}(z,z') \;=\; m\,{\rm L}_{m-1,n}(z,z'),\quad \frac{\partial}{\partial z'}\,{\rm L}_{m,n}(z,z') \;=\; n\,{\rm L}_{m,n-1}(z,z'),\label{difflag2d}
\end{eqnarray}
and, furthermore, the Laguerre 2D polynomials satisfy the following recurrence relations
\begin{eqnarray}
{\rm L}_{m+1,n}(z,z') &=& z\,{\rm L}_{m,n}(z,z')-n\,{\rm L}_{m,n-1}(z,z'),\nonumber\\ {\rm L}_{m,n+1}(z,z') &=& z'\,{\rm L}_{m,n}(z,z')-m\,{\rm L}_{m-1,n}(z,z'),
\end{eqnarray}
as was derived in [7--9] and as can be easily seen from (\ref{lag2d}) or (\ref{lag2d1}).

The Laguerre 2D polynomials (\ref{lag2d1}) are related to the generalized Laguerre (or Laguerre-Sonin)\footnote{N.Ya. 'Sonin' is often
written in the French form N.J. 'Sonine' under which this Russian mathematician of 19-th to 20-th
century became known in Western Europe.} polynomials ${\rm L}_n^{\nu}(u)$ by
\begin{eqnarray}
{\rm L}_{m,n}(z,z') &=& (-1)^n n! z^{m-n}{\rm L}_n^{m-n}(zz')\nonumber\\ &=& (-1)^m m!
z'^{\,n-m}{\rm L}_m^{n-m}(zz'),\label{lag2d2}
\end{eqnarray}
that explains the given name. In most physical applications the second complex variable
$z'$ is complex conjugated to the first variable $z$ that means $z'=z^*$ but for generality we leave open this specialization and consider $z$ and $z'$ as two independent complex (or sometimes
real) variables.

The operators $z-\frac{\partial}{\partial z'}$ and $z'-\frac{\partial}{\partial z}$ which play a role in (\ref{lag2d}) are commutative that means
\begin{eqnarray}
\left[z-\frac{\partial}{\partial z'},z'-\frac{\partial}{\partial z}\right] &\equiv & \left(z-\frac{\partial}{\partial z'}\right)\left(z'-\frac{\partial}{\partial z}\right)-\left(z'-\frac{\partial}{\partial z}\right)\left(z-\frac{\partial}{\partial z'}\right) \;=\;0,
\end{eqnarray}
and their powers can be disentangled (all multiplication operators stand in front of the differential operators) that using the explicit form of the Laguerre 2D polynomials (\ref{lag2d1}) leads to the following operator identity
\begin{eqnarray}
\left(z-\frac{\partial}{\partial z'}\right)^m \left(z'-\frac{\partial}{\partial z}\right)^n &=&
\sum_{k=0}^m\sum_{l=0}^n\frac{(-1)^{k+l} m!n!}{k!(m-k)!\,l!(n-l)!}\,{\rm L}_{m-k,n-l}(z,z')\frac{\partial^{k+l}}{\partial z'^k\partial z^l}.\label{laguerreop}
\end{eqnarray}
\textcolor{red}{It} is applicable to arbitrary functions $f(z,z')$ and provides then functional identities such as (\ref{lag2d}) in application to the function $f(z,z')=1$.
For its derivation we used in addition the reordering relation of differentiation and multiplication operators (\ref{laguerreop}) (e.g., \cite{w12}, Eq. (A.6) there)
\begin{eqnarray}
\frac{\partial^l}{\partial z^l}z^k &=& \sum_{j=0}^{\{k,l\}}\frac{k!l!}{j!(k-j)!(l-j)!}z^{k-j}\frac{\partial^{l-j}}{\partial z^{l-j}},
\end{eqnarray}
which is well known in quantum optics (transition from antinormal to normal ordering of boson creation and annihilation operators) and can be proved by complete induction but relation (\ref{laguerreop}) can be also directly proved by complete induction.

The (special) Laguerre 2D polynomials ${\rm L}_{m,n}(z,z')\equiv {\rm L}_{m,n}\left({\sf
I};z,z'\right)$ are the special case ${\sf U}={\sf I}$ of the (general) Laguerre 2D
polynomials ${\rm L}_{m,n}\left({\sf U};z,z'\right)$ where ${\sf U}$ is a general 2D
matrix and ${\sf I}$ denotes the 2D unit matrix 
[10--12]. Together with the
general Hermite 2D polynomials ${\rm H}_{m,n}\left({\sf U};x,y\right)$, the general
Laguerre 2D polynomials ${\rm L}_{m,n}\left({\sf U};z,z'\right)$ form a unified object
which can be transformed from one to the other form by a special unitary matrix ${\sf Z}$
which transforms the real coordinates $(x,y)$ to the pair of complex coordinates
$(z\equiv x+\I y,z'= z^*\equiv x-\I y)$. It seems to be not an overestimation to say that
the appearance of the generalized Laguerre polynomials ${\rm L}_n^{\nu}(u)$ in
applications most often in the form of ${\rm L}_{m,n}(z,z^*)$ leads to the conclusion
that the usual generalized Laguerre polynomials ${\rm L}_n^{\nu}(u)$ are the radial
rudiments of the Laguerre 2D polynomials ${\rm L}_{m,n}\left(r\E^{\I\varphi},
r\E^{-\I\varphi}\right)$ in polar coordinates $\left(z\equiv r\E^{\I\varphi}, z^*\equiv
r\E^{-\I\varphi}\right)$ with $u\equiv r^2$. Their orthonormalization on the positive
semi-axis $0\le u \equiv r^2 < +\infty$ with weight proportional to $\E^{-u}$ supports
this conclusion.

The Laguerre 2D polynomials are related to products of Hermite polynomials by \cite{w4}
(the special case $m=n$ is given in \cite{bate2,soni} but with an error by an absent factor $2^{2n}$ on the right-hand side)
\begin{eqnarray}
{\rm L}_{m,n}\left(x+\I y,x-\I y\right) &=& (-1)^n n!r^{m-n}{\rm
L}_{n}^{m-n}\left(r^2\right)\E^{\I(m-n)\varphi} \nonumber\\ &=&
(-1)^n\sum_{j=0}^{m+n}\left(\frac{\I}{2}\right)^{m+n-j}\,{\rm P}_{j}^{(m-j,n-j)}(0)\,{\rm
H}_j(x){\rm H}_{m+n-j}(y),\label{lagjacobiherm}
\end{eqnarray}
and the inversion is
\begin{eqnarray}
{\rm H}_m(x){\rm H}_n(y) &=& \I^n \sum_{j=0}^{m+n}2^j\, {\rm P}_{j}^{(m-j,n-j)}(0){\rm
L}_{j,m+n-j}\left(x+\I y,x-\I y\right) \nonumber\\ &=& \I^n \sum_{j=0}^{m+n}(-2)^j\, {\rm
P}_{j}^{(m-j,n-j)}(0)j!r^{m+n-2j}{\rm L}_j^{m+n-2j}\left(r^2\right)\E^{\I(2j-m-n)\varphi},
\end{eqnarray}
where the coefficients are essentially given by the Jacobi polynomials ${\rm
P}_j^{(\alpha,\beta)}(u)$ for argument $u=0$. Background of these formulae is the
relation [8--11]
\begin{eqnarray}
z^mz^{*n} \;=\;\left(x+\I y\right)^m\left(x-\I y\right)^n \;=\; (-1)^n\sum_{j=0}^{m+n}\left(\frac{\I}{2}\right)^{m+n-j}\,{\rm P}_{j}^{(m-j,n-j)}(0)(2x)^j (2y)^{m+n-j}, \label{zzstarxy}
\end{eqnarray}
and its inversion
\begin{eqnarray}
(2x)^m (2y)^n \;=\; \left(z+z^*\right)^m\left(-\I\left(z-z^*\right) \right)^n
\;=\; \I^n\sum_{j=0}^{m+n}2^j\,{\rm P}_{j}^{(m-j,n-j)}(0)z^{j}z^{*\,m+n-j},
\end{eqnarray}
and the application of the integral operator $\exp\left(-\frac{\partial^2}{\partial z \partial z^*}\right)=\exp\left(-\frac{1}{4}\left(\frac{\partial^2}{\partial x^2}+\frac{\partial^2}{\partial y^2}\right)\right)$ to them (see also Section 2). In special case $m=n$ the last two formulae make the transition to
\begin{eqnarray}
\left(zz^*\right)^n &=& \left(x^2+y^2\right)^n \;=\;\sum_{k=0}^{n}\frac{n!}{k!(n-k)!}x^{2k}y^{2(n-k)} \nonumber\\ &=& \sum_{k=0}^{n}(-1)^k2^{2k}\,{\rm P}_{2k}^{(n-2k,n-2k)}(0)x^{2k}y^{2(n-k)},
\end{eqnarray}
and to
\begin{eqnarray}
(4xy)^n &=& \left(-\I\left(z^2-z^{*2}\right)\right)^n \;=\; \I^n\sum_{k=0}^{n}\frac{(-1)^k n!}{k!(n-k)!}z^{2k}z^{*2(n-k)} \nonumber\\ &=& \I^n\sum_{k=0}^{n}2^{2k}\,{\rm P}_{2k}^{(n-2k,n-2k)}(0)z^{2k}z^{*2(n-k)},\label{xyzzstar}
\end{eqnarray}
from which by comparison of the different representations follows
\begin{eqnarray}
{\rm P}_{2k}^{(n-2k,n-2k)}(0)= \frac{(-1)^k n!}{2^{2k}k!(n-k)!}, \quad {\rm P}_{2k+1}^{(n-2k-1,n-2k-1)}(0)=0,
\end{eqnarray}
for Jacobi polynomials with equal upper indices specialized by argument $u=0$. The Jacobi polynomials with equal upper indices are also called Ultraspherical polynomials ${\rm P}_n^{(\alpha,\alpha)}(u)$ and are related to Gegenbauer polynomials ${\rm C}_n^{\nu}(u)$ by
\begin{eqnarray}
{\rm P}_n^{(\alpha,\alpha)}(u) \;\equiv\; \frac{(2\alpha)!(n+\alpha)!}{\alpha!(n+2\alpha)!}{\rm
C}_n^{\alpha+\frac{1}{2}}(u), \quad {\rm C}_n^{\nu}(u) \;\equiv\;  \frac{\left(\nu-\frac{1}{2}\right)!(n+2\nu-1)!} {(2\nu-1)!\left(n+\nu-\frac{1}{2}\right)!}\,{\rm P}_n^{\left(\nu-\frac{1}{2},\nu-\frac{1}{2}\right)}(u).
\end{eqnarray}
Clearly, formulae (\ref{zzstarxy})--(\ref{xyzzstar}) remain to be true if $z^*\rightarrow z'$ is not the complex conjugate to $z$ in which case $x$ and $y$ become complex numbers.

Different methods of derivation of generating functions are presented in the monographs
\cite{rain,sriv}. The problem of determination of the basic generating function for simple
Laguerre 2D and Hermite 2D polynomials was solved in 
[9--12, 18--20]. A more difficult problem is the determination of generating functions for products of two Laguerre 2D polynomials or of a Laguerre 2D and a Hermite 2D polynomial. In \cite{w6}, we derived some special generating functions for products of two Laguerre 2D polynomials. The corresponding generating functions with general 2D matrices ${\sf U}$ as parameters in these polynomials are fairly complicated \cite{w6}. In present paper, we derive by an operational approach the generating functions for products of two special Laguerre 2D polynomials, for products of two Hermite 2D polynomials and for the mixed case of a product of a Laguerre 2D with a Hermite 2D polynomial (also called bilinear generating functions). This corresponds to the formula of Mehler (e.g., \cite{bate2}, 10.13 (22) and below in Section 3) which is the bilinear generating function for the product of two usual Hermite polynomials. We begin in next Section with a short representation of the analogical 1D case of Hermite polynomials and discuss in Section 3 their bilinear generating function and continue in Sections 4--7 with the corresponding derivations for the Laguerre 2D and Hermite 2D cases. In Section 8 we derive a summation formula over Laguerre 2D polynomials which can be considered as intermediate step to the mentioned generating functions but possesses also its own importance in applications. Sections 9 and 10 are concerned with the further illumination of two factorizations of two different bilinear generating functions.

The operators which play a role in one of the definitions of Hermite and Laguerre (1D and
2D) polynomials are Gaussian convolutions and possess a relation to the Lie group
$SU(1,1)$. Using operator disentanglement for $SU(1,1)$ we may derive operational
relations which provide a useful tool for the derivation of the considered generating
functions. This is presented in Appendix A.

\setcounter{chapter}{2}
\setcounter{equation}{0}
\section*{2. Hermite polynomials and their alternative definition as a 1D analogue to Laguerre 2D polynomials}

Hermite polynomials ${\rm H}_n(x)$\footnote{In Russian literature, Hermite polynomials are often called Chebyshev-Hermite polynomials.} (e.g., 
[1--4]) can be defined in analogy to Laguerre 2D polynomials (\ref{lag2d}), at least, in two well-suited equivalent ways by
\begin{eqnarray}
{\rm H}_n(x) & \equiv & \exp\left(-\frac{1}{4}\frac{\partial^2}{\partial
x^2}\right)(2x)^n \underbrace{\exp\left(\frac{1}{4}\frac{\partial^2}{\partial
x^2}\right)\, 1}_{=1} \nonumber\\ & = & (-1)^n\exp\left(x^2\right)\frac{\partial^n}{\partial
x^n}\exp\left(-x^2\right)\, 1 \nonumber\\ &=& \left(2x-\frac{\partial}{\partial x}\right)^n\,1,\label{hermite0}
\end{eqnarray}
that leads to the following alternative definition from which results immediately the explicit representation
\begin{eqnarray}
{\rm H}_n(x) &=& \exp\left(-\frac{1}{4}\frac{\partial^2}{\partial x^2}\right)(2x)^n \;=\;
\sum_{l=0}^{\left[\frac{n}{2}\right]}\frac{(-1)^l n!}{l!(n-2l)!}(2x)^{n-2l},\label{hermite}
\end{eqnarray}
with the inversion
\begin{eqnarray}
(2x)^n &=& \exp\left(\frac{1}{4}\frac{\partial^2}{\partial x^2}\right)\,{\rm H}_n(x) \;=\;
\sum_{l=0}^{\left[\frac{n}{2}\right]}\frac{n!}{l!(n-2l)!}\,{\rm H}_{n-2l}(x).\label{herminv}
\end{eqnarray}
The definition (\ref{hermite}) which is little known (see \cite{bell}, 5.3, pp. 159/160) and which was occasionally used in older time is an alternative one to the well-known Rodrigues-type definition given in second line in (\ref{hermite0}) and it has found new attention and its fixed place in literature only in recent time 
[27--33]. It rests here on the operator identities $ \left(n=0,1,2,\ldots\right)$
\begin{eqnarray}
\left(2x-\frac{\partial}{\partial x}\right)^n &=& \left(\exp\left(-\frac{1}{4}\frac{\partial^2}{\partial
x^2}\right)2x\,\exp\left(\frac{1}{4}\frac{\partial^2}{\partial
x^2}\right)\right)^n \nonumber\\ &=& \left(-\exp\left(x^2\right)\frac{\partial}{\partial x}\exp\left(-x^2\right)\right)^n \nonumber\\ &=& \sum_{l=0}^n \frac{n!}{l!(n-l)!}\,{\rm H}_{n-l}(x)\left(-\frac{\partial}{\partial x}\right)^l,\label{hermiteop}
\end{eqnarray}
which can be applied to arbitrary functions $f(x)$ of a (in general, complex) variable $x$ and is applied in (\ref{hermite0}) to the function $f(x)=1$. Similar considerations were made for the operators in the definition (\ref{lag2d}) of Laguerre 2D polynomials where we have as background of the two alternatives the identity of the operators in (\ref{lag2d}) leading to the operator identity (\ref{laguerreop}).

Such alternative definitions of a sequence of polynomials $p_n(x),\,(n=0,1,2,\ldots)$ as in (\ref{hermite}) are possible in every case if the generating function possesses a special form and vice versa as follows
\begin{eqnarray}
\sum_{n=0}^\infty \frac{t^n}{n!}p_n(x) = \E^{axt}f(t),\quad \Leftrightarrow \quad p_n(x) = f\left(\frac{1}{a}\frac{\partial}{\partial x}\right)(ax)^n,\label{expgf}
\end{eqnarray}
with an arbitrary function $f(t)$. The proof is relatively simple and is here omitted. In case of the Hermite polynomials one has $f(t)=\exp\left(-t^2\right)$ with the parameter $a=2$. Other examples are binomials, higher-order Hermite polynomials (Gould-Hopper polynomials \cite{sriv} with little applications up to now), Bernoulli polynomials and Euler polynomials the last related to Hyperbolic Secant function. The analogous alternative definition of Laguerre 2D polynomials is given in first line of (\ref{lag2d1}) in comparison to a more conventional one in second line. For some proofs the alternative definitions possess advantages but we will not and cannot state this generally.

Generalized Laguerre polynomials ${\rm L}_n^{\nu}(u)$ form a peculiar case with respect to the generating functions since according to (\ref{lag2d2}) they are properly rudiments of Laguerre 2D polynomials where they are involved in the form ${\rm L}_n^{m-n}(u)$. Indeed, in this combination of indices they possess a known generating function of the form (\ref{expgf}) (e.g., 
[1--4, 26]. Some other kinds of generating functions for ${\rm L}_n^{\nu}(u)$ with fixed $\nu$ can be also found in cited literature with the possibility of an operational definition as follows ($a=-1$ and $f(t)=(1+t)^m$ in (\ref{expgf}))
\begin{eqnarray}
&& \sum_{n=0}^\infty t^n\,{\rm L}_n^{m-n}(u)=\sum_{n=0}^\infty \frac{t^n}{n!}\,n!\,{\rm L}_n^{m-n}(u) =\E^{-ut}\left(1+t\right)^{m},\label{opdeflag}
\end{eqnarray}
from which follows (compare with (\ref{lag2d1}))
\begin{eqnarray}
n!\,{\rm L}_n^{m-n}(u) &=& \left(1-\frac{\partial}{\partial u}\right)^{m}(-u)^n \;=\; \sum_{j=0}^{\{m,n\}}\frac{m!n!}{j!(m-j)!(n-j)!}(-u)^{n-j}, \nonumber\\ n!\,{\rm L}_n^{\nu}(u) &=& \left(1-\frac{\partial}{\partial u}\right)^{\nu +n}(-u)^n \;=\; n!\,\left(1-\frac{\partial}{\partial u}\right)^{\nu} \underbrace{{\rm L}_n^0(u)}_{\equiv\, {\rm L}_n(u)}, \label{opdeflag1}
\end{eqnarray}
where $\left(1-\frac{\partial}{\partial u}\right)^{\nu}$ for arbitrary real values of $\nu$ is defined by their Taylor series in powers of $\frac{\partial}{\partial u}$. The known expansion of ${\rm L}_n^{m-n}(u)$ in powers of $u$ can be immediately checked from the operational definition. The yet remaining classical polynomials are the Jacobi polynomials ${\rm P}_n^{(\alpha,\beta)}(z)$ with their special cases. They do not possess an operational definition according to the scheme mentioned here for Hermite and Laguerre polynomials. However, a more complicated form of an operational definition was found in \cite{w9} (Appendix A there) for the special case of ultraspherical polynomials ${\rm P}_n^{(\alpha,\alpha)}(z)$ and thus also for their equivalent Gegenbauer polynomials and their special cases of Legendre and of Chebyshev polynomials. It is unclear up to now whether or not and in which form there exists an extension to the general case of Jacobi polynomials ${\rm P}_n^{(\alpha,\beta)}(z)$ with $\alpha \neq \beta$.

Formulae (\ref{expgf}) are closely related to the so-called umbral calculus \cite{roman} in its simplest form and it seems that an essential part of this symbolic calculus rests on the duality of the linear functionals of the delta function $\delta(x)$ and its derivatives $\delta^{(m)}(x)$ to the monomials $x^n$ according to $\left((-1)^m\delta^{(m)}(x),\frac{x^n}{n!}\right) \equiv \frac{(-1)^m}{n!}\int_{-\infty}^{+\infty}dx\,\delta^{(m)}(x)x^n =\delta_{m,n}$ and to the introduction of a symbolic notation for the calculation with the corresponding algebra which, however, does not bring a great relief in comparison to the direct calculation with these functionals.

\setcounter{chapter}{3}
\setcounter{equation}{0}
\section*{3. Generating function for products of two Hermite polynomials (Mehler formula)}

Besides the well-known generating function for Hermite polynomials
\begin{eqnarray}
\sum_{n=0}^{\infty}\frac{t^n}{n!}{\rm H}_n(x) &=& \exp\left(2tx-t^2\right),\label{genherm}
\end{eqnarray}
in applications, in particular, in quantum optics of the harmonic oscillator the following bilinear
generating function for products of two Hermite polynomials with equal indices but
different arguments plays an important role (formula of Mehler; see, e.g. \cite{bate2},
10.13 (22))
\begin{eqnarray}
\sum_{n=0}^{\infty}\frac{t^n}{2^n n!}{\rm H}_n(x){\rm H}_n(y) &=& \frac{1}{\sqrt{1-t^2}}
\exp\left(\frac{2txy-t^2\left(x^2+y^2\right)}{1-t^2}\right)
\nonumber\\ &=& \frac{1}{\sqrt{1+t}}\exp\left(\frac{t\left(x+y\right)^2}{2(1+t)}\right) \frac{1}{\sqrt{1-t}}\exp\left(-\frac{t\left(x-y\right)^2}{2(1-t)}\right). \quad (|t|< 1).\quad \label{mehler}
\end{eqnarray}
We represented here the right-hand side additionally in a sometimes useful
factorization. This factorization is connected with the following identity (see
\cite{prud2}, 4.5.2 (5), p. 641)
\begin{eqnarray}
{\rm H}_n(x){\rm H}_n(y) &=& \frac{1}{2^{n}}\sum_{k=0}^n\frac{(-1)^k n!} {k!(n-k)!}\,{\rm
H}_{2(n-k)}\left(\frac{x+y}{\sqrt{2}}\right)\,{\rm H}_{2k}\left(\frac{x-y}{\sqrt{2}}\right) \nonumber\\ &=& 2^{n}n!\sum_{k=0}^n(-1)^{n-k}\,{\rm L}_{n-k}^{-\frac{1}{2}}\left(\frac{(x+y)^2}{2}\right)\,{\rm L}_{k}^{-\frac{1}{2}} \left(\frac{(x-y)^2}{2} \right), \label{factherm}
\end{eqnarray}
and thus with a coordinate transformation. In the special case $y=0$ of (\ref{mehler}) using
\begin{eqnarray}
{\rm H}_{2m}(0)=\frac{(-1)^m(2m)!}{m!},\quad {\rm H}_{2m+1}(0),\quad \left(m=0,1,2,\ldots\right), \label{h2m0}
\end{eqnarray}
it provides (see, e.g., \cite{datt2}, Eq. (75) presented there, however, with a less usual definition of Hermite polynomials)
\begin{eqnarray}
\sum_{m=0}^\infty \frac{(-1)^m}{m!}\left(\frac{t}{2}\right)^{2m}\,{\rm H}_{2m}(x)
&=& \frac{1}{\sqrt{1-t^2}}\exp\left(-\frac{(tx)^2}{1-t^2}\right), \label{genfevenhermpol}
\end{eqnarray}
which is a generating function for even Hermite polynomials and by differentiation with respect to variable $x$ using $\frac{\partial}{\partial x}\,{\rm H}_{n}(x)=2n\,{\rm H}_{n-1}(x)$
\begin{eqnarray}
\sum_{m=0}^\infty \frac{(-1)^m}{m!}\left(\frac{t}{2}\right)^{2m+1}\,{\rm H}_{2m+1}(x)
&=& \frac{tx}{\big(\sqrt{1-t^2}\,\big)^{3}}\exp\left(-\frac{(tx)^2}{1-t^2}\right),\label{genfoddhermpol}
\end{eqnarray}
the corresponding generating functions for odd Hermite polynomials. We mention here that both generating functions (\ref{genfevenhermpol}) and (\ref{genfoddhermpol}) are not contained in the otherwise very comprehensive and impressive work \cite{bate2} but instead of these are two other ones for even and odd Hermite polynomials which can easily be obtained by separating the even and odd part in the most well-known usual generating function for Hermite polynomials (\ref{genherm}). Let us give them since a small mistake is there ($\sqrt{2}$ on the right-hand sides in \cite{bate2} (10.13, (20) and (21)) has to be changed to $2$)
\begin{eqnarray}
\sum_{m=0}^\infty\frac{t^{2m}}{(2m)!}\,{\rm H}_{2m}(x) &=& \sum_{n=0}^\infty\frac{\frac{1}{2}\left(t^n+(-t)^n\right)}{n!}\,{\rm H}_n(x) \;=\; \exp\left(-t^2\right)\,{\rm ch}(2tx), \nonumber\\ \sum_{m=0}^\infty\frac{t^{2m+1}}{(2m+1)!}\,{\rm H}_{2m+1}(x) &=& \sum_{n=0}^\infty\frac{\frac{1}{2}\left(t^n-(-t)^n\right)}{n!}\,{\rm H}_n(x) \;=\; \exp\left(-t^2\right)\,{\rm sh}(2tx).
\end{eqnarray}
Correct formulae of this kind one may find in \cite{soni} (chap. V, p. 252).

In the limiting transition $t\rightarrow +1$ setting $t=1-\frac{\vep}{2}$ in (\ref{mehler}), one obtains for $\vep\rightarrow 0$ the following completeness relation for the Hermite polynomials
\begin{eqnarray}
\sum_{n=0}^{\infty}\frac{1}{2^n n!}{\rm H}_n(x){\rm H}_n(y) &=&
\exp\left(x^2+y^2-xy\right)\lim_{\vep\rightarrow
+0}\frac{1}{\sqrt{\vep}}\exp\left(-\frac{(x-y)^2}{\vep}\right)\nonumber\\ &=&
\sqrt{\pi}\exp\left(\frac{x^2+y^2}{2}\right)\delta(x-y).
\end{eqnarray}
As it is known this indicates a way for the introduction of Hermite functions ${\rm h}_n(x)$
as follows
\begin{eqnarray}
{\rm h}_n(x) &\equiv &
\frac{1}{\pi^{\frac{1}{4}}}\exp\left(-\frac{x^2}{2}\right)\frac{{\rm
H}_n(x)}{\sqrt{2^nn!}},
\end{eqnarray}
which are complete and orthonormalized according to
\begin{eqnarray}
&& \sum_{n=0}^\infty {\rm h}_n(x){\rm h}_n(y) = \delta(x-y),\quad
\int_{-\infty}^{+\infty}dx\,{\rm h}_{m}(x){\rm h}_{n}(x)=\delta_{m,n}.
\end{eqnarray}
This underlines the great importance of the bilinear generating function (\ref{mehler}).
A proof of (\ref{mehler}) can be given, for example, using an operational formula derived
in Appendix A (Eq. (\ref{op2}) there) or using (\ref{factherm}) in connection with the
generating function for generalized Laguerre polynomials that shifts it to the proof
of (\ref{factherm}) (the more direct proof of (\ref{mehler}) in \cite{rain} occupies one full page on pp. 197/198).

\setcounter{chapter}{4}
\setcounter{equation}{0}
\section*{4. Generating functions for Laguerre 2D polynomials}

In the derivation of generating functions for Laguerre 2D polynomials ${\rm
L}_{m,n}(z,z')$ it does not play a role whether $z'$ is complex conjugated to $z$ or not. The
Laguerre 2D polynomials ${\rm L}_{m,n}(z,z')$ possess the following symmetry and scaling
property
\begin{eqnarray}
&&{\rm L}_{m,n}(z,z') = {\rm L}_{n,m}(z',z),\quad {\rm L}_{m,n}\left(\kappa
z,\frac{z'}{\kappa}\right)=\kappa^{m-n}{\rm L}_{m,n}(z,z'),\label{scal}
\end{eqnarray}
for arbitrary complex $\kappa$. Thus the scaling transformation of the variables $z\rightarrow \kappa z,z'\rightarrow \frac{z'}{\kappa}$ preserves the Laguerre 2D polynomials up to factors
$\kappa^{m-n}$ which reduces to pure phase factors $\E^{\I(m-n)\chi}$ for $|\kappa|=1$ or $\kappa=\E^{\I\chi}$ with special case $\chi=\pi$ or $\kappa=-1$ and therefore ${\rm L}_{m,n}(-z,-z')=(-1)^{m-n}{\rm L}(z,z')$. In case that $z'$ is the complex conjugated variable $z'=z^{*}$ to $z$ and therefore not independent on $z$ we may only consider the transformation $z\rightarrow \E^{\I\chi}z,z^{*}\rightarrow \E^{-\I\chi}z^{*}$ where the Laguerre 2D polynomials are preserved up to phase factors $\E^{\I(m-n)\chi}$. The scaling transformations can be used to make a certain check of the final results for generating functions.

The following generating function for Laguerre 2D polynomials is easily to obtain from the alternative definition of ${\rm L}_{m,n}(z,z')$ in (\ref{lag2d1}) in the following way
\begin{eqnarray}
\sum_{m=0}^\infty\sum_{n=0}^\infty \frac{s^mt^n}{m!n!}\,{\rm L}_{m,n}(z,z') &=& \exp\left(-\frac{\partial^2}{\partial z \partial z'}\right) \sum_{m=0}^\infty\sum_{n=0}^\infty \frac{s^mt^n}{m!n!}\,z^m z'^n \nonumber\\ &=& \exp\left(-\frac{\partial^2}{\partial z \partial z'}\right) \exp\left(sz+tz'\right) \nonumber\\ &=&
\exp\left(sz+tz'-st\right),\label{simplgenflag2d}
\end{eqnarray}
in particular for $s=z,t=z'$ and for $s=z',t=z$
\begin{eqnarray}
\sum_{m=0}^\infty\sum_{n=0}^\infty \frac{z^m z'^{\,n}}{m!n!}\,{\rm L}_{m,n}(z,z') &=& \exp\left(z^2+z'^{2}-zz'\right), \nonumber\\
\sum_{m=0}^\infty\sum_{n=0}^\infty \frac{z'^{\,m} z^n}{m!n!}\,{\rm L}_{m,n}(z,z') &=& \exp\left(zz'\right).
\end{eqnarray}

As an intermediate step to the generating
(\ref{simplgenflag2d}) one may consider the generating function with summation over only one of the indices in the special Laguerre 2D polynomials
\begin{eqnarray}
\sum_{n=0}^\infty \frac{t^n}{n!}\,{\rm L}_{m,n}(z,z') &=& \left(z-t\right)^m\exp\left(tz'\right),
\end{eqnarray}
which is to obtain analogously to (\ref{simplgenflag2d}) using the first definition of these polynomials in (\ref{lag2d}). Written by the usual generalized Laguerre polynomials according to (\ref{lag2d2}) this provides
\begin{eqnarray}
\sum_{n=0}^{\infty}\left(-\frac{z}{t}\right)^{n}\,{\rm L}_n^{m-n}(zz') &=& (-1)^m \frac{m!}{(zz')^{\,m}} \sum_{n=0}^{\infty}\frac{(t z')^n}{n!}\,{\rm L}_m^{n-m}(zz')
\nonumber\\ &=& \left(1-\frac{t}{z}\right)^m \exp\left(\frac{t}{z}zz'\right),
\end{eqnarray}
which after substitutions $\frac{t}{z} \rightarrow -t'$ and $zz'\rightarrow u$ becomes identical with the known generating function (\ref{opdeflag}) for usual Laguerre polynomials \cite{bate2}.

\setcounter{chapter}{5}
\setcounter{equation}{0}
\section*{5. Generating functions for products of Laguerre 2D polynomials}

We now calculate the basic generating function for the product of two special Laguerre 2D
polynomials. Using the first of the definitions in (\ref{lag2d}) which corresponds to
the alternative definition of the Hermite polynomials in first line of (\ref{hermite}),
we quickly proceed as follows (remind (\ref{scal}))
\begin{eqnarray}
&& \sum_{m=0}^\infty\sum_{n=0}^\infty \frac{s^mt^n}{m!n!}{\rm L}_{m,n}(z,z'){\rm
L}_{m,n}(w',w)\nonumber\\ &=& \exp\left(-\frac{\partial^2}{\partial z\partial
z'}-\frac{\partial^2}{\partial w\partial
w'}\right)\sum_{m=0}^\infty\frac{(szw')^m}{m!}\sum_{n=0}^\infty\frac{(twz')^n}{n!}
\nonumber\\ &=& \exp\left(-\frac{\partial^2}{\partial z\partial
z'}-\frac{\partial^2}{\partial w\partial w'}\right)\exp\left(szw'+twz'\right)\nonumber\\
&=& \exp\left(-\frac{\partial^2}{\partial z\partial
z'}\right)\exp\left(szw'+twz'-stzz'\right)\nonumber\\ &=&
\exp\left(ww'\right)\exp\left(-\frac{\partial^2}{\partial z\partial z'}\right)
\exp\left\{-st\left(z-\frac{w}{s}\right)\left(z'-\frac{w'}{t}\right)\right\},\label{genllmn}
\end{eqnarray}
where in third step \textcolor{red}{it} is used that $\exp\left(szw'+twz'\right)$ is an eigenstate of the operator $\frac{\partial^2}{\partial w\partial w'}$ to eigenvalues $stzz'$. The remaining
step to solve is essentially the convolution $\exp\left(-\frac{\partial^2}{\partial
z\partial z'}\right)\exp\left(-stzz'\right)$ with succeeding displacement $z\rightarrow
z-\frac{w}{s},z'\rightarrow z'-\frac{w'}{t}$ in the result of this convolution. We can
calculate this convolution via two-dimensional Fourier transformation but in Appendix A
(Eq. (\ref{op4})) we derive an operational formula which allows a more direct approach
and which is very useful for similar calculations. With the result of this convolution we
obtain
\begin{eqnarray}
\sum_{m=0}^\infty\sum_{n=0}^\infty \frac{s^mt^n}{m!n!}{\rm L}_{m,n}(z,z'){\rm L}_{n,m}(w,w') &=& \sum_{m=0}^\infty\sum_{n=0}^\infty \frac{s^mt^n}{m!n!}{\rm L}_{m,n}(z,z'){\rm L}_{m,n}(w',w)  \nonumber\\ &=& \frac{1}{1-st}\exp\left(\frac{szw'+twz'-st(zz'+ww')}{1-st}\right). \label{gflag2dpr}
\end{eqnarray}
This result for the generating function can be factorized in the following form
\begin{eqnarray}
\frac{1}{1-st}\exp\left(\frac{szw'+twz'-st(zz'+ww')}{1-st}\right) &=&
\frac{1}{1+\sqrt{st}}\exp\left(\frac{\left(\sqrt{s}\,z+\sqrt{t}\,w\right)
\left(\sqrt{t}\,z'+\sqrt{s}\,w'\right)}{2\left(1+\sqrt{st}\right)}\right)\nonumber\\ && \cdot
\frac{1}{1-\sqrt{st}}\exp\left(-\frac{\left(\sqrt{s}\,z-\sqrt{t}\,w\right)
\left(\sqrt{t}\,z'-\sqrt{s}\,w'\right)}{2\left(1-\sqrt{st}\right)}\right),\nonumber\\
\label{gflag2dprfact}
\end{eqnarray}
in analogy to the formula of Mehler (\ref{mehler}) and a fully analogous derivation of this formula by coordinate transformations is possible. In Section 10 we derive such a decomposition of the product of two Laguerre 2D polynomials with the same indices but different arguments which provides a further insight into the factorization in (\ref{gflag2dpr}) according to (\ref{gflag2dprfact}). If we substitute in (\ref{gflag2dpr}) $s\rightarrow \frac{s}{w'},t\rightarrow \frac{t}{w}$
and if we use in the obtained modified generating function the limiting transition
\begin{eqnarray}
\lim_{{\sscr |w|\rightarrow \infty, |w'|\rightarrow \infty}}\frac{\,{\rm L}_{m,n}(w',w)}{w'^{m}w^n} =1,
\end{eqnarray}
then we find by this limiting procedure from (\ref{gflag2dpr}) the generating function
(\ref{simplgenflag2d}).

Expressed by usual generalized Laguerre polynomials, relation (\ref{gflag2dpr}) takes on
the following forms
\begin{eqnarray}
&&\sum_{m=0}^\infty\sum_{n=0}^\infty \frac{n!}{m!}s^mt^n\left(zw'\right)^{m-n}{\rm
L}_n^{m-n}(zz'){\rm L}_n^{m-n}(ww') \nonumber\\ &=& \sum_{m=0}^\infty\sum_{n=0}^\infty
s^mt^n\left(-\frac{z'}{w'}\right)^{n-m}{\rm L}_m^{n-m}(zz'){\rm L}_n^{m-n}(ww')\nonumber\\
&=& \frac{1}{1-st}\exp\left(\frac{szw'+twz'-st(zz'+ww')} {1-st}\right),
\end{eqnarray}
where on the left-hand side relations (\ref{lag2d2}) were used.

Similar to the case of Hermite polynomials, we now make the limiting transition
$s\rightarrow 1,t\rightarrow 1$ in the generating function (\ref{gflag2dpr}) and make for
this purpose the specializations of complex conjugation $z'=z^*,\;w'=w^*$ of the variables.
With $s=1-\frac{\vep}{2},t=1-\frac{\vep}{2}$, we obtain
\begin{eqnarray}
&&\sum_{m=0}^\infty\sum_{n=0}^\infty \frac{1}{m!n!}{\rm L}_{m,n}(z,z^*)\left({\rm
L}_{m,n}(w,w^*)\right)^* \nonumber\\ &=&
\exp\left(zz^*+ww^*-\frac{1}{2}\left(zw^*+wz^*\right)\right)\lim_{\vep\rightarrow +1}
\frac{1}{\vep}\exp\left(-\frac{\left(z-w\right)\left(z^*-w^*\right)}{\vep}\right)
\nonumber\\ &=& \pi \exp\left(\frac{zz^*+ww^*}{2}\right)\delta\left(z-w,z^*-w^*\right),
\end{eqnarray}
where $\delta(z,z^*)=\delta(x+\I y,x-\I y)\equiv \delta(x)\delta(y)$
denotes the two-dimensional delta function in representation by a pair of complex conjugated variables. This relation (and also the orthogonality relations) suggest to introduce the following Laguerre 2D functions ${\rm l}_{m,n}(z,z^*)$ by \cite{w2}
\begin{eqnarray}
{\rm l}_{m,n}(z,z^*) &\equiv &
\frac{1}{\sqrt{\pi}}\exp\left(-\frac{zz^*}{2}\right)\frac{{\rm
L}_{m,n}(z,z^*)}{\sqrt{m!n!}},\label{compll2d}
\end{eqnarray}
which are complete and orthonormalized according to
\begin{eqnarray}
&&\sum_{m=0}^\infty\sum_{n=0}^\infty {\rm l}_{m,n}(z,z^*)\left({\rm
l}_{m,n}(w,w^*)\right)^* = \delta(z-w,z^*-w^*),\nonumber\\ && \int \frac{\I}{2} dz\wedge
dz^* \left({\rm l}_{k,l}(z,z^*)\right)^*{\rm l}_{m,n}(z,z^*) = \delta_{k,m}\delta_{l,n}.
\label{orthl2d}
\end{eqnarray}
Herein, $\frac{\I}{2}dz\wedge dz^* = dx\wedge dy$ is the area element of the complex
plane and the integration goes over the whole complex plane. Relations (\ref{compll2d})
and (\ref{orthl2d}) can be used for the expansion of functions of two variables in Laguerre 2D
polynomials or Laguerre 2D functions.

By forming derivatives of (\ref{gflag2dpr}) with respect to the variables, one can derive
related formulae. Furthermore, by specialization of the variables, we obtain new
formulae. For example, due to ${\rm L}_{m,n}(0,0)=(-1)^n n!\delta_{m,n}$ (see (\ref{lag2d1z0}))
we obtain from (\ref{gflag2dpr}) setting $w=w'=0$ and substituting $ st\rightarrow t$
\begin{eqnarray}
\sum_{n=0}^\infty \frac{(-t)^n}{n!}{\rm L}_{n,n}(z,z') &=& \sum_{n=0}^\infty
t^n{\rm L}_n(zz') \;=\; \frac{1}{1-t}\exp\left(-\frac{tzz'}{1-t}\right).\label{genlagpol}
\end{eqnarray}
This is the well-known generating function for the usual Laguerre polynomials ${\rm
L}_n(u)\equiv {\rm L}_n^0(u)$. Using the definition (\ref{lag2d1}) this can be calculated also
by
\begin{eqnarray}
\sum_{n=0}^\infty \frac{(-t)^n}{n!}{\rm L}_{n,n}(z,z') &=& \exp\left(-\frac{\partial^2}{\partial z \partial z'}\right) \sum_{n=0}^\infty \frac{(-tzz')^n}{n!} \nonumber\\ &=& \exp\left(-\frac{\partial^2}{\partial z \partial z'}\right) \exp\left(-tzz'\right), \label{genlagpol1}
\end{eqnarray}
that is the convolution of two 2D Gaussian functions (see Appendix, formulae (\ref{op4}) or (\ref{convtwogauss}) with corresponding substitutions) which provides the result on the
right-hand side of (\ref{genlagpol}).

\setcounter{chapter}{6}
\setcounter{equation}{0}
\section*{6. Generating functions for products of two Hermite 2D polynomials}

Special Hermite 2D polynomials ${\rm H}_{m,n}(x,y)$ are the special case ${\sf U}={\sf
I}$ of general Hermite 2D polynomials ${\rm H}_{m,n}({\sf U};x,y)$ and are products of
usual Hermite polynomials according to
\begin{eqnarray}
{\rm H}_{m,n}(x,y) \equiv {\rm H}_{m,n}({\sf I};x,y) = {\rm H}_m(x){\rm H}_n(y).
\label{specherm2d}
\end{eqnarray}
The generating function for products of two special Hermite 2D polynomials defined by (\ref{specherm2d}) factorizes
\begin{eqnarray}
\sum_{m=0}^\infty \sum_{n=0}^\infty \frac{s^m t^n}{2^{m+n}m!n!}{\rm H}_{m,n}(x,y){\rm
H}_{m,n}(u,v) &=& \sum_{m=0}^\infty \frac{s^m}{2^{m} m!}{\rm H}_{m}(x){\rm
H}_{m}(u)\sum_{n=0}^\infty \frac{t^n}{2^{n}n!}{\rm H}_{n}(y){\rm H}_n(v)\nonumber\\ &=&
\sum_{m=0}^\infty \frac{s^m}{2^{m} m!}{\rm H}_{m,m}(x,u)\,\sum_{n=0}^\infty \frac{t^n}{2^{n}n!}{\rm H}_{n,n}(y,v),\label{gfh2dpr0}
\end{eqnarray}
and is easily to obtain by using explicitly the Mehler formula (\ref{mehler}) with the result
\begin{eqnarray}
&& \sum_{m=0}^\infty \sum_{n=0}^\infty \frac{s^m t^n}{2^{m+n}m!n!}{\rm H}_{m,n}(x,y){\rm
H}_{m,n}(u,v) \nonumber\\ &=&
\frac{1}{\sqrt{\left(1-s^2\right)\left(1-t^2\right)}}\exp\left(\frac{2sxu-s^2\left(
x^2+u^2\right)}{1-s^2}+\frac{2tyv-t^2\left(y^2+v^2\right)}{1-t^2}\right).\label{gfh2dpr}
\end{eqnarray}
In special case $u=v=0$ using ${\rm H}_{2k}(0)=\frac{(-1)^k(2k)!}{k!},\,{\rm H}_{2k+1}(0) = 0,\,(k=0,1,\ldots)$ we find from this
\begin{eqnarray}
&& \sum_{k=0}^\infty\frac{(-1)^k}{k!}\left(\frac{s}{2}\right)^{2k}\,{\rm H}_{2k}(x)
\sum_{l=0}^\infty\frac{(-1)^l}{l!}\left(\frac{t}{2}\right)^{2n}\,{\rm H}_{2l}(y) \nonumber\\ &=&
\frac{1}{\sqrt{\left(1-s^2\right)\left(1-t^2\right)}} \exp\left(-\frac{(sx)^2}{1-s^2}-
\frac{(ty)^2}{1-t^2}\right),
\end{eqnarray}
that is the product of two generating functions for even Hermite polynomials of the form
(\ref{genfevenhermpol}). By differentiations of this formula with respect to variables $x$
or $y$ one finds also the cases for odd Hermite polynomials and the mixed case of even and odd
Hermite polynomials but due to factorization this belongs already to the primary stage of the generating functions (\ref{genfevenhermpol}) and (\ref{genfoddhermpol}).

\setcounter{chapter}{7}
\setcounter{equation}{0}
\section*{7. Generating function for products of Laguerre 2D and Hermite 2D polynomials and
for Laguerre 2D polynomials with even indices as its special case}

Sometimes in quantum-optical calculations, one has to evaluate the following double sum
of the form of a generating function for products of special Laguerre 2D with special
Hermite 2D polynomials the last defined by (\ref{specherm2d}) and can quickly proceed to the following stage
\begin{eqnarray}
&&\sum_{m=0}^\infty \sum_{n=0}^\infty \frac{s^m t^n}{\sqrt{2^{m+n}}\,m!n!}{\rm
L}_{m,n}(z,z'){\rm H}_{m,n}(u,v)\nonumber\\ &=& \exp\left\{-\frac{\partial^2}{\partial
z\partial z'}-\frac{1}{4}\left(\frac{\partial^2}{\partial u^2}+\frac{\partial^2}{\partial
v^2}\right)\right\}\sum_{m=0}^\infty \sum_{n=0}^\infty \frac{s^m
t^n}{\sqrt{2^{m+n}}\,m!n!}z^mz'^n(2u)^m(2v)^n \nonumber\\ &=&
\exp\left\{-\frac{\partial^2}{\partial z\partial
z'}-\frac{1}{4}\left(\frac{\partial^2}{\partial u^2}+\frac{\partial^2}{\partial
v^2}\right)\right\}\exp\left(\sqrt{2}\left(szu+tz'v\right)\right).\label{gflagh1}
\end{eqnarray}
Here we have two possibilities to continue the calculations. Using two times the
generating function for Hermite polynomials, we obtain from (\ref{gflagh1})
\begin{eqnarray}
&&\sum_{m=0}^\infty \sum_{n=0}^\infty \frac{s^m t^n}{\sqrt{2^{m+n}}\,m!n!}{\rm
L}_{m,n}(z,z'){\rm H}_{m,n}(u,v)\nonumber\\ &=& \exp\left(u^2+v^2\right)
\exp\left(-\frac{\partial^2}{\partial z\partial
z'}\right)\exp\left\{-\frac{s^2}{2}\left(z-\frac{\sqrt{2}\,u}{s}\right)^2-
\frac{t^2}{2}\left(z'-\frac{\sqrt{2}\,v}{t}\right)^2\right\},\label{gflagh2}
\end{eqnarray}
and using first the generating function for Laguerre 2D polynomials, we find
\begin{eqnarray}
&&\sum_{m=0}^\infty \sum_{n=0}^\infty \frac{s^m t^n}{\sqrt{2^{m+n}}\,m!n!}{\rm
L}_{m,n}(z,z'){\rm H}_{m,n}(u,v) \nonumber\\ &=& \exp\left(zz'\right)
\exp\left\{-\frac{1}{4}\left(\frac{\partial^2}{\partial u^2}+\frac{\partial^2}{\partial
v^2}\right)\right\}\exp\left\{-2st\left(u-\frac{z'}{\sqrt{2}\,s}\right)
\left(v-\frac{z}{\sqrt{2}\,t}\right)\right\}.\label{gflagh3}
\end{eqnarray}
The remaining problem is to calculate the two-dimensional convolutions in (\ref{gflagh2})
or (\ref{gflagh3}). Both convolutions can be accomplished using auxiliary formulae
prepared in Appendix A. The result is
\begin{eqnarray}
&&\sum_{m=0}^\infty \sum_{n=0}^\infty \frac{s^m t^n}{\sqrt{2^{m+n}}\,m!n!}{\rm
L}_{m,n}(z,z'){\rm H}_{m,n}(u,v) = \frac{1}{\sqrt{1-s^2t^2}} \nonumber\\ && \quad \cdot
\exp\left(\frac{2\sqrt{2}\left(suz+tvz'
+st\left(svz+tuz'\right)\right)-s^2z^2-t^2z'^2-4stuv-2s^2t^2\left(zz'
+u^2+v^2\right)}{2\left(1-s^2t^2\right)}\right). \nonumber\\ \label{gflagh4}
\end{eqnarray}
The complexity of this generating function finds a simple explanation to which we say some words
at the end of this Section.

In the special case $z=z'=0$ using (\ref{lag2d1z0}), we obtain from (\ref{gflagh4}) with the substitution $st\rightarrow -t$ the Mehler formula (\ref{mehler}). In the special
case $u=v=0$ using ${\rm H}_{2k}(0)=\frac{(-1)^k(2k)!}{k!}, {\rm H}_{2k+1}(0)=0,$ we find from (\ref{gflagh4})
\begin{eqnarray}
\sum_{k=0}^\infty \sum_{l=0}^\infty \frac{(-1)^{k+l}s^{2k} t^{2l}}{k!l!2^{k+l}}{\rm
L}_{2k,2l}(z,z') &=& \frac{1}{\sqrt{1-s^2t^2}}\exp\left(-\frac{s^2z^2+t^2z'^2+2s^2t^2zz'}
{2\left(1-s^2t^2\right)}\right)\nonumber\\ && \hspace{-30mm} \;=\; \frac{1}{\sqrt{1-st}}\exp\left( -\frac{\left(sz+tz'\right)^2}{4(1-st)}\right)\frac{1}{\sqrt{1+st}}\exp\left(
-\frac{\left(sz-tz'\right)^2}{4(1+st)}\right), \label{gflagh5}
\end{eqnarray}
where, in addition, a factorization is given. Using the generating function (\ref{genfevenhermpol}) for even Hermite polynomials we get from this factorization the identity
\begin{eqnarray}
\sum_{k=0}^\infty \sum_{l=0}^\infty \frac{(-1)^{k+l}s^{2k} t^{2l}}{k!l!2^{k+l}}{\rm
L}_{2k,2l}(z,z') &=& \sum_{m=0}^\infty \frac{(-1)^m}{m!}\left(\frac{\sqrt{st}}{2}\right)^{2m}\,{\rm H}_{2m}\left(\frac{1}{2}\left(\sqrt{\frac{s}{t}}z+\sqrt{\frac{t}{s}}z'\right)\right) \nonumber\\
&& \cdot \sum_{n=0}^\infty \frac{(-1)^n}{n!}\left(\I\frac{\sqrt{st}}{2}\right)^{2n}\,{\rm H}_{2n}\left(-\frac{\I}{2}\left(\sqrt{\frac{s}{t}}z-\sqrt{\frac{t}{s}}z'\right)\right), \qquad
\end{eqnarray}
which we consider from another point of view in Section 9 providing thus some better understanding for it.

For $s=0$ or for $t=0$ we obtain from (\ref{gflagh4}) the common generating function for Hermite polynomials with obvious substitutions.
On the other side, the generating function (\ref{gflagh5}) for Laguerre 2D polynomials with even indices can be more directly obtained from (we substitute here $s^2\rightarrow -\sigma,t^2\rightarrow -\tau$)
\begin{eqnarray}
&&\sum_{k=0}^\infty \sum_{l=0}^\infty \frac{\sigma^k}{2^k k!}\frac{\tau^l}{2^l l!}\,{\rm L}_{2k,2l}\left(z,z'\right) \;=\;
\exp\left(-\frac{\partial^2}{\partial z \partial z'}\right)\sum_{k=0}^\infty \frac{1}{k!} \left(\frac{\sigma}{2}z^2\right)^k\sum_{l=0}^\infty \frac{1}{l!}\left(\frac{\tau}{2}z'^2\right)^l \nonumber\\ &=& \exp\left(-\frac{\partial^2}{\partial z \partial z'}\right)
\exp\left(\frac{\sigma}{2}z^2+\frac{\tau}{2}z'^2\right) \;=\; \exp\left(\frac{\sigma}{2}\left(z-\frac{\partial}{\partial z'}\right)^2\right)\exp\left(-\frac{\partial^2}{\partial z \partial z'}\right)\exp\left(\frac{\tau}{2}z'^2\right) \nonumber\\ &=& \exp\left(\frac{\sigma}{2}z^2-\sigma z\frac{\partial}{\partial z'}+\frac{\sigma}{2}\frac{\partial^2}{\partial z'^2}\right)\exp\left(\frac{\tau}{2}z'^2\right),
\nonumber\\ &=& \exp\left(\frac{\sigma}{2}z^2-\sigma z\frac{\partial}{\partial z'}\right)\frac{1}{\sqrt{1-\sigma \tau}}\,\exp\left(\frac{\tau z'^2}{2\left(1-\sigma \tau\right)}\right),
\end{eqnarray}
where we used the identity (\ref{op1hermpol}) in the Appendix in special case $n=0$.
Accomplishing the last operation of argument displacement of variable $z'$ by the operator $\exp\left(-\sigma z\frac{\partial}{\partial z'}\right)$ we obtain the generating function
\begin{eqnarray}
\sum_{k=0}^\infty \sum_{l=0}^\infty \frac{\sigma^k\tau^l}{2^{k+l} k! l!}\,{\rm L}_{2k,2l}\left(z,z'\right) &=&\frac{1}{\sqrt{1-\sigma\tau}}\,\exp\left(\frac{\sigma z^2+\tau z'^2-2\sigma\tau zz'}{2\left(1-\sigma\tau\right)}\right), \end{eqnarray}
which with obvious substitutions ($s^2\rightarrow -\sigma,t^2\rightarrow -\tau$) is identical with (\ref{gflagh5}). By differentiation of this generating function with respect to variables $z$ and (or) $z'$ one obtains generating functions for Laguerre 2D polynomials with odd (or even and odd) indices.

We mention yet that expressed by the usual generalized Laguerre polynomials according to (\ref{lag2d2}) and by applying the doubling formula for the argument of the Gamma function the left-hand side of (\ref{gflagh5}) can be written ($\left(-\frac{1}{2}\right)! = \sqrt{\pi}\,$)
\begin{eqnarray}
\sum_{k=0}^\infty \sum_{l=0}^\infty \frac{(-1)^{k+l}s^{2k} t^{2l}}{k!l!2^{k+l}}{\rm
L}_{2k,2l}(z,z') &=& \sum_{k=0}^\infty \sum_{l=0}^\infty \frac{\left(l-\frac{1}{2}\right)!s^{2k}t^{2l}}{\left(-\frac{1}{2}\right)!k!}
\left(-\frac{z^2}{2}\right)^{k-l}\,{\rm L}_{2l}^{2(k-l)}(zz') \nonumber\\ &=& \sum_{k=0}^\infty \sum_{l=0}^\infty \frac{\left(k-\frac{1}{2}\right)!s^{2k}t^{2l}}{\left(-\frac{1}{2}\right)!l!}
\left(-\frac{z'^2}{2}\right)^{l-k}\,{\rm L}_{2k}^{2(l-k)}(zz'),
\end{eqnarray}
the last by symmetry of the Laguerre 2D polynomials or using (\ref{lag2d2}).

As already mentioned in the Introduction the special Laguerre 2D and Hermite 2D polynomials can be combined in one whole object of polynomials ${\rm L}_{m,n}({\sf U},z,z')$ or ${\rm H}_{m,n}({\sf V};x,y)$, alternatively, with general 2D matrices ${\sf U}$ and ${\sf V}$ and can be transformed into each other in this form where a special matrix ${\sf Z}$ plays a main role 
[8--11]. The generating function (\ref{gflagh4}) belongs to a special case where expressed by the more general Hermite 2D or Laguerre 2D polynomials the polynomials ${\rm L}_{m,n}({\sf I};z,z')$ and ${\rm L}_{m,n}({\sf Z};u,v)$ or ${\rm H}_{m,n}({\sf Z}^{-1};z,z')$ and ${\rm H}_{m,n}({\sf I};u,v)$ (${\sf I}$ is unit matrix) are joined in one formula and such cases become complicated written in components of the matrices. Therefore formula (\ref{gflagh4}) may also play a role as nontrivial special case of more general generating functions for arbitrary different matrices ${\sf U}$ and ${\sf V}$ in the polynomials and we have checked (\ref{gflagh4}) also numerically.

\setcounter{chapter}{8}
\setcounter{equation}{0}
\section*{8. A set of simple sums over products of
Laguerre 2D polynomials}

We now consider a set of simple (in the sense of not double!) sums over products of Laguerre 2D polynomials with two free indices $(m, n=0,1,2,\ldots)$ as follows
\begin{eqnarray}
&&\sum_{k=0}^\infty \frac{(-t)^k}{k!}{\rm L}_{m,k}(z,z'){\rm L}_{k,n}(w,w') \nonumber\\
&=& \exp\left(-\frac{\partial^2}{\partial z \partial z'}-\frac{\partial^2}{\partial w
\partial w'}\right)\sum_{k=0}^\infty \frac{(-twz')^k}{k!}z^mw'^{\,n}\nonumber\\ &=&
\exp\left(-\frac{\partial^2}{\partial z \partial z'}-\frac{\partial^2}{\partial w
\partial w'}\right)\exp\left(-twz'\right)z^m w'^{\,n}\nonumber\\ &=&
\exp\left\{-t\left(z'-\frac{\partial}{\partial
z}\right)\left(w-\frac{\partial}{\partial w'}\right)\right\}z^mw'^{\,n}\nonumber\\ &=&
\exp\left(-twz'\right)\exp\left\{t\left(w\frac{\partial}{\partial
z}+z'\frac{\partial}{\partial w'}\right)\right\}\exp\left(-t\frac{\partial^2}{\partial
z\partial w'}\right)z^mw'^{\,n}\nonumber\\ &=&
\exp\left(-twz'\right)\exp\left\{t\left(w\frac{\partial}{\partial
z}+z'\frac{\partial}{\partial w'}\right)\right\}\left(\sqrt{t}\,\right)^{m+n}{\rm
L}_{m,n}\left(\frac{z}{\sqrt{t}}, \frac{w'}{\sqrt{t}}\right).\label{lagsum0}
\end{eqnarray}
The last two steps of making the argument displacements and using the scaling property
(\ref{scal}) lead to the following final representations (among other possible ones)
\begin{eqnarray}
\sum_{k=0}^\infty \frac{(-t)^k}{k!}{\rm L}_{m,k}(z,z'){\rm L}_{k,n}(w,w') &=& \sum_{k=0}^\infty \frac{(-t)^k}{k!}{\rm L}_{m,k}(z,z'){\rm L}_{n,k}(w',w) \nonumber\\ &=& \exp\left(-twz'\right)\left(\sqrt{t}\, \right)^{m+n}{\rm L}_{m,n}\left(\frac{z+tw}{\sqrt{t}}, \frac{w'+t z'}{\sqrt{t}}\right)\nonumber\\
&=& \exp\left(-twz'\right)t^n{\rm L}_{m,n}\left(z+tw,z'+\frac{w'}{t}\right).
\label{lagsum1}
\end{eqnarray}
Expressed by generalized Laguerre polynomials using (\ref{lag2d2}), this relation takes
on the form (compare a similar form in \cite{prud2}, (chap. 5.11.5. Eq.(2.)))
\begin{eqnarray}
&& z^mw'^{\,n}\sum_{k=0}^\infty k!\left(-\frac{t}{zw'}\right)^k{\rm L}_k^{m-k}(zz') {\rm
L}_k^{n-k}(ww')\nonumber\\ &=& \frac{(-1)^{m+n}m!n!}{z'^{\,m}w^n}\sum_{k=0}^\infty
\frac{(-twz')^k}{k!}{\rm L}_m^{k-m}(zz'){\rm L}_n^{k-n}(ww')\nonumber\\
&=& \exp\left(-twz'\right)n!(-t)^n\left(z+tw\right)^{m-n}{\rm
L}_n^{m-n}\left(\frac{(z+tw)(w'+tz')}{t}\right) \nonumber\\ &=&
\exp\left(-twz'\right)m!(-t)^m\left(w'+tz'\right)^{n-m}{\rm
L}_m^{n-m}\left(\frac{(z+tw)(w'+tz')}{t}\right).\label{lagsum2}
\end{eqnarray}
In the special case $(w,w')=(z,z')$ one finds from (\ref{lagsum1})
\begin{eqnarray}
\sum_{k=0}^\infty \frac{(-t)^k}{k!}{\rm L}_{m,k}(z,z'){\rm L}_{k,n}(z,z') &=& \exp\left(-tzz'\right)\left(\sqrt{t}\, \right)^{m+n}{\rm L}_{m,n}\left(\frac{1+t}{\sqrt{t}}\,z, \frac{1+t}{\sqrt{t}}\,z'\right),\label{lagsum3}
\end{eqnarray}
or expressed by the generalized Laguerre polynomials after division of (\ref{lagsum2}) by $z^mw'^n$
\begin{eqnarray}
&& \sum_{k=0}^\infty k!\left(-\frac{t}{u}\right)^k{\rm L}_k^{m-k}(u) {\rm
L}_k^{n-k}(u)\nonumber\\ &=& \frac{(-1)^{m+n}m!n!}{u^{m+n}}\sum_{k=0}^\infty
\frac{(-tu)^k}{k!}{\rm L}_m^{k-m}(u){\rm L}_n^{k-n}(u)\nonumber\\
&=& \exp\left(-tu\right)n!\left(-\frac{t}{u}\right)^n\left(1+t\right)^{m-n}{\rm
L}_n^{m-n}\left(\frac{(1+t)^2}{t}u\right) \nonumber\\ &=&
\exp\left(-tu\right)m!\left(-\frac{t}{u}\right)^m\left(1+t\right)^{n-m}{\rm
L}_m^{n-m}\left(\frac{(1+t)^2}{t}u\right),\label{lagsum4}
\end{eqnarray}
where we made the substitution $u \equiv zz'$.

In most representations of orthogonal polynomials, one can only find generating functions
for products of generalized Laguerre polynomials where the upper indices are parameters
and are not involved in the summations (e.g., \cite{bate2} (chap. 10.12.(20)) and \cite{rain,sriv}). Using the limiting
relation
\begin{eqnarray}
\lim_{\vep\rightarrow 0}\vep^k{\rm
L}_k^{n-k}\left(\frac{u}{\vep}\right)=\frac{(-1)^k}{k!}u^k,
\end{eqnarray}
and substituting $w'=\frac{w''}{\vep}$ in (\ref{lagsum2}), we obtain by limiting
procedure $\vep\rightarrow 0$
\begin{eqnarray}
\frac{1}{z^m}\sum_{k=0}^\infty \frac{(-tw)^k}{k!}{\rm L}_{m,k}(z,z') &\equiv &
\sum_{k=0}^\infty\left(\frac{tw}{z}\right)^k{\rm L}_{k}^{m-k}\left(zz'\right) \nonumber\\ &=&
\exp\left(-\frac{tw}{z}zz'\right)\left(1+\frac{tw}{z}\right)^m,
\end{eqnarray}
where $w''$ disappeared. This is with substitutions the known relation (\ref{opdeflag}) which in \cite{bate2} as mentioned is classified under generating functions (chap. 10.12, Eq. (19)) and which in our representation by special Laguerre 2D polynomials proves to be one of their basic generating functions with simple summation.

The sums over products of Laguerre 2D polynomials (\ref{lagsum1}) or (\ref{lagsum2}) possess proper importance for sum evaluations which sometimes arise when working with these polynomials. They also form a partial result on the way to the evaluation of the generating function (\ref{gflag2dpr}) for the product of two Laguerre 2D polynomials. Taking in (\ref{lagsum1}) the special case $m=n$ and multiplying it by $\frac{s^m}{m!}$ and forming then the sum over $m$ we obtain with substitution $t\rightarrow -t$ and using the well-known generating function for Laguerre polynomials ${\rm L}_m(u)\equiv {\rm L}_m^{0}(u)$
\begin{eqnarray}
&& \sum_{m=0}^\infty\sum_{n=0}^\infty \frac{s^mt^n}{m!n!}{\rm L}_{m,n}(z,z'){\rm
L}_{m,n}(w',w)\nonumber\\ &=& \exp\left(twz'\right)\sum_{m=0}^\infty (st)^m {\rm
L}_m\left(-\frac{(z-tw)(w'-tz')}{t}\right)\nonumber\\ &=& \exp\left(twz'\right)
\frac{1}{1-st}\exp\left(\frac{s\left(z-tw\right)\left(w'-tz'\right)}{1-st}\right).
\label{gflag2d}
\end{eqnarray}
Joining herein the two exponential functions we see that the right-hand side of (\ref{gflag2d}) is equal to the right-hand side of (\ref{gflag2dpr}) as it is necessary and thus we have calculated here this generating function in a second way.

\setcounter{chapter}{9}
\setcounter{equation}{0}
\section*{9. Factorization of generating function for simple Laguerre 2D polynomials with even indices}
We illuminate now a cause for the possible factorization in the generating function (\ref{gflagh5}) for Laguerre 2D polynomials with even indices. For this purpose we make in (\ref{gflagh5}) the substitutions
\begin{eqnarray}
&& z= \sqrt{\frac{t}{s}} \left(x+\I y\right), \quad z'= \sqrt{\frac{s}{t}} \left(x-\I y\right), \quad x=\frac{1}{2}\left(\sqrt{\frac{s}{t}}z+\sqrt{\frac{t}{s}}z'\right), \quad y = -\frac{\I}{2}\left(\sqrt{\frac{s}{t}}z-\sqrt{\frac{t}{s}}z'\right),
\nonumber\\&& \frac{\partial}{\partial z} =\frac{1}{2}\sqrt{\frac{s}{t}}\left(\frac{\partial}{\partial x} - \I\frac{\partial}{\partial y}\right), \quad \frac{\partial}{\partial z'} =\frac{1}{2}\sqrt{\frac{t}{s}}\left(\frac{\partial}{\partial x} + \I\frac{\partial}{\partial y}\right), \quad \frac{\partial^2}{\partial z \partial z'} = \frac{1}{4}\left(\frac{\partial^2}{\partial x^2}+\frac{\partial^2}{\partial y^2}\right).\label{transzzprxy}
\end{eqnarray}
Therefore
\begin{eqnarray}
&& \sum_{k=0}^\infty \sum_{l=0}^\infty \frac{(-1)^{k+l}s^{2k} t^{2l}}{k!l!2^{k+l}}{\rm
L}_{2k,2l}(z,z') \;=\; \exp\left(-\frac{1}{4}\left(\frac{\partial^2}{\partial x^2}+\frac{\partial^2}{\partial y^2}\right)\right)\sum_{k=0}^\infty \sum_{l=0}^\infty\frac{(-1)^{k+l}}{k!l!}\left(\frac{st}{2}\right)^{k+l} \nonumber\\ && \hspace{60mm} \cdot\frac{1}{2}\left((x+\I y)^{2k}(x-\I y)^{2l}+(x-\I y)^{2k}(x+\I y)^{2l}\right)
\nonumber\\ &=& \exp\left(-\frac{1}{4}\left(\frac{\partial^2}{\partial x^2}+\frac{\partial^2}{\partial y^2}\right)\right)\sum_{k=0}^\infty \sum_{l=0}^\infty \frac{(-1)^{k+l}}{k!l!}\left(\frac{st}{2}\right)^{k+l}\sum_{n=0}^{2(k+l)}2^{2n}{\rm P}_{2n}^{(2(k-n),2(l-n)}(0)x^{2(k+l-n)}(\I y)^{2n} \nonumber\\ &=& \exp\left(-\frac{1}{4}\left(\frac{\partial^2}{\partial x^2}+\frac{\partial^2}{\partial y^2}\right)\right)\sum_{m=0}^\infty\sum_{n=0}^\infty (-1)^{m} \left(\frac{st}{2}\right)^{m+n}x^{2m}y^{2n} \sum_{k=0}^{m+n}\frac{2^{2n}{\rm P}_{2n}^{(2(k-n),2(m-k))}(0)}{k!(m+n-k)!} \nonumber\\ &=&  \exp\left(-\frac{1}{4}\left(\frac{\partial^2}{\partial x^2}+\frac{\partial^2}{\partial y^2}\right)\right)\sum_{m=0}^\infty\sum_{n=0}^\infty \frac{(-1)^{m}}{m!n!}(st)^{m+n}x^{2m}y^{2m}
\nonumber\\ &=& \exp\left(-\frac{1}{4}\frac{\partial^2}{\partial x^2}\right)\exp\left(-st x^2\right)\exp\left(-\frac{1}{4}\frac{\partial^2}{\partial y^2}\right)\exp\left(st y^2\right),
\end{eqnarray}
that can be written
\begin{eqnarray}
\sum_{k=0}^\infty \sum_{l=0}^\infty \frac{(-1)^{k+l}s^{2k} t^{2l}}{k!l!2^{k+l}}{\rm
L}_{2k,2l}(z,z') &=& \sum_{m=0}^\infty\frac{1}{m!}\left(-\frac{\sqrt{st}}{2}\right)^{2m}
{\rm H}_{2m}(x)\sum_{n=0}^\infty\frac{1}{n!}\left(\frac{\sqrt{st}}{2}\right)^{2n}\,{\rm H}_{2n}(y) \nonumber\\ &=& \frac{1}{\sqrt{1-st}}\exp\left(-\frac{st x^2}{1-st}\right)\frac{1}{\sqrt{1+st}}\exp\left(\frac{st y^2}{1+st}\right),\qquad
\end{eqnarray}
where the alternative definition of Hermite polynomials in (\ref{hermite0}) is applied and
where an apparently unknown sequence of finite sum evaluations
\begin{eqnarray}
\sum_{k=0}^{m+n}\frac{2^{2n}{\rm P}_{2n}^{(2(k-n),2(m-k))}(0)}{k!(m+n-k)!} &=& \frac{1}{(2m)!(2n)!} \sum_{k=0}^{m+n}\frac{(2k)!(2(m+n-k))!}{k!(m+n-k)!}2^{2k}{\rm P}_{2k}^{(2(m-k),2(n-k))}(0) \nonumber\\ &=& \frac{2^{m+n}}{m!n!},\qquad \left(m=0,1,\ldots;\,n=0,1,\ldots\,\right),\label{sumid}
\end{eqnarray}
is inserted. These sum identities are proved already by the obvious equivalence (\ref{gflagh5}) and are easily to check for small $(m,n)$. A direct independent proof we did not make but probably it is possible by complete induction.
If we go back to the variables $(z,z')$ according to (\ref{transzzprxy}) we have the factorization (\ref{gflagh5}).

\setcounter{chapter}{10}
\setcounter{equation}{0}
\section*{10. Identities for products of two Laguerre 2D polynomials with different arguments}

We derive here a decomposition of the generating function for an entangled product of two equal Laguerre 2D polynomials with different arguments into a product of two simple generating function for Laguerre 2D polynomials which in its further application leads immediately to the given factorization in the generating function (\ref{gflag2dpr}) and provides a certain explanation for it. The considerations are in some sense similar to the considerations in the previous Section.

We now make in (\ref{gflag2dpr}) the following substitutions of variables
\begin{eqnarray}
&& z=\left(\frac{t}{s}\right)^{\frac{1}{4}}\frac{x+y}{\sqrt{2}\,},\quad w =\left(\frac{s}{t}\right)^{\frac{1}{4}}\frac{x-y}{\sqrt{2}\,},\quad  z'= \left(\frac{s}{t}\right)^{\frac{1}{4}}\frac{x'+y'}{\sqrt{2}\,},\quad w' =\left(\frac{t}{s}\right)^{\frac{1}{4}}\frac{x'-y'}{\sqrt{2}\,},\qquad \label{zwzwxyxy}
\end{eqnarray}
with the inversion
\begin{eqnarray}
&& x = \frac{\left(\frac{s}{t}\right)^{\frac{1}{4}}z+\left(\frac{t}{s}\right)^{\frac{1}{4}}w} {\sqrt{2}},\quad x' = \frac{\left(\frac{t}{s}\right)^{\frac{1}{4}}z' +\left(\frac{s}{t}\right)^{\frac{1}{4}}w'}{\sqrt{2}},\nonumber\\&&
y = \frac{\left(\frac{s}{t}\right)^{\frac{1}{4}}z-\left(\frac{t}{s} \right)^{\frac{1}{4}}w}{\sqrt{2}},\quad y' = \frac{\left(\frac{t}{s}\right)^{\frac{1}{4}}z'-\left(\frac{s}{t}\right)^{\frac{1}{4}}w'}{\sqrt{2}},
\label{xyxyzwzw}
\end{eqnarray}
from which follows for the operators of differentiation
\begin{eqnarray}
&& \frac{\partial}{\partial z} =\frac{1}{\sqrt{2}}\,\left(\frac{s}{t}\right)^{\frac{1}{4}}\left(\frac{\partial}{\partial x}+\frac{\partial}{\partial y}\right),\quad \frac{\partial}{\partial z'} =\frac{1}{\sqrt{2}}\,\left(\frac{t}{s}\right)^{\frac{1}{4}}\left(\frac{\partial}{\partial x'}+\frac{\partial}{\partial y'}\right),\nonumber\\ && \frac{\partial}{\partial w} =\frac{1}{\sqrt{2}}\,\left(\frac{t}{s}\right)^{\frac{1}{4}}\left(\frac{\partial}{\partial x}-\frac{\partial}{\partial y}\right),\quad \frac{\partial}{\partial w'} =\frac{1}{\sqrt{2}}\,\left(\frac{s}{t}\right)^{\frac{1}{4}}\left(\frac{\partial}{\partial x'}-\frac{\partial}{\partial y'}\right).
\end{eqnarray}
As a consequence we find in transformed coordinates (see also (\ref{genllmn})) with application of formula (\ref{convtwogauss}) in the Appendix for the evaluation of the Gaussian integrals
\begin{eqnarray}
&=& \sum_{m=0}^\infty\sum_{n=0}^\infty \frac{s^mt^n}{m!n!}\,{\rm L}_{m,n}(z,z')\,{\rm L}_{n,m}(w,w') \nonumber\\ &=& \exp\left(-\frac{\partial^2}{\partial x \partial x'}-\frac{\partial^2}{\partial y \partial y'}\right)\sum_{m=0}^\infty \frac{\left(\sqrt{st}(x+y)(x'-y')\right)^m}{2^mm!}\sum_{n=0}^\infty \frac{\left(\sqrt{st}(x-y)(x'+y')\right)^n}{2^nn!} \nonumber\\ &=& 
\exp\left(-\frac{\partial^2}{\partial x \partial x'}\right)\exp\left(\sqrt{st}\,xx'\right)
\exp\left(-\frac{\partial^2}{\partial y \partial y'}\right)\exp\left(-\sqrt{st}\,yy'\right),
\end{eqnarray}
that using the generating function (\ref{genlagpol}) for simple Laguerre 2D polynomials with equal indices can be written
\begin{eqnarray}
&&\sum_{m=0}^\infty\sum_{n=0}^\infty \frac{s^mt^n}{m!n!}\,{\rm L}_{m,n}(z,z')\,{\rm L}_{n,m}(w,w') \nonumber\\ &=& \sum_{m=0}^\infty\frac{\left(\sqrt{st}\,\right)^m}{m!}\,{\rm L}_{m,m}(x,x')\,\sum_{n=0}^\infty \frac{\left(-\sqrt{st}\,\right)^n}{n!}\,{\rm L}_{n,n}(y,y') \nonumber\\ &=& \frac{1}{1+\sqrt{st}}\exp\left(\frac{\sqrt{st}\,xx'}{1+\sqrt{st}}\right)\frac{1}
{1-\sqrt{st}}\exp\left(-\frac{\sqrt{st}\,yy'}{1-\sqrt{st}}\right).\qquad
\end{eqnarray}
If we go back on the right-hand side to the primary variables $(z,z',w,w')$ according to (\ref{xyxyzwzw}) we arrive at the given factorization (\ref{gflag2dprfact}). One may look at this as to an alternative derivation of the bilinear generating function (\ref{gflag2dpr}).

\setcounter{chapter}{11}
\setcounter{equation}{0}
\section*{11. Comparison of the two alternative definitions in the derivation of the generating function for Hermite polynomials}
It is not possible to say generally which of the two alternative definitions of Hermite polynomials in Section 2 and of Laguerre 2D polynomials in Section 1 are better to work with. This depends on the problem and sometimes also on the taste and the former experience of the user. We demonstrate this in the simplest case of the derivation of the generating function (\ref{genherm}) for Hermite polynomials. Using definition in the first line of (\ref{hermite0}) we calculate
\begin{eqnarray}
\sum_{n=0}^{\infty}\frac{t^n}{n!}{\rm H}_n(x) &=& \exp\left(-\frac{1}{4}\frac{\partial^2}{\partial x^2}\right)\sum_{n=0}^{\infty}\frac{t^n}{n!}(2x)^n \;=\; \exp\left(-\frac{1}{4}\frac{\partial^2}{\partial x^2}\right)\exp\left(2tx\right) \nonumber\\ &=& \exp\left(2tx-t^2\right),
\end{eqnarray}
and using definition in the second line of (\ref{hermite0})
\begin{eqnarray}
\sum_{n=0}^{\infty}\frac{t^n}{n!}{\rm H}_n(x) &=& \exp\left(x^2\right)\sum_{n=0}^{\infty}\frac{t^n}{n!}(-1)^n\frac{\partial^n}{\partial x^n}\exp\left(-x^2\right) \;=\; \exp\left(x^2\right)\exp\left(-t\frac{\partial}{\partial x}\right)\exp\left(-x^2\right) \nonumber\\ &=& \exp\left(x^2\right)\exp\left(-\left(x-t\right)^2\right).
\end{eqnarray}
In the first case we use that $\exp(2tx)$ is eigenfunction of the operator $\frac{1}{4}\frac{\partial^2}{\partial x^2}$ to the eigenvalue $t^2$ and in the second case that
$\exp\left(-t\frac{\partial}{\partial x}\right)$ is the displacement operator of the argument $x$ of a function $f(x)$ to $f(x-t)$. Here both derivations are equally simple. However, in case of the inversion formula (\ref{herminv}) and of formula (\ref{lagjacobiherm}), for example, the alternative definition (\ref{hermite}) seems to be more suited.

The alternative definitions of Hermite (1D and 2D) and Laguerre 2D polynomials extends the arsenal of possible approaches to problems of their application and one should have for disposal the new method in the same way as the former methods.

\setcounter{chapter}{12}
\setcounter{equation}{0}
\section*{12. Conclusion}

We have derived and discussed generating functions for the product of two special Laguerre 2D or Hermite 2D polynomials and for the mixed case of such products. In our derivations we preferred the first (operational) definition of the Laguerre 2D polynomials (\ref{lag2d1}) from the two alternative ones given in (\ref{lag2d}) which as it seems to us is advantageous for this purpose. This is due to the separation of the same operator $\exp\left(-\frac{\partial^2}{\partial z \partial z'}\right)$ applied to $z^mz'^{\,n}$ for all polynomials $\,{\rm L}_{m,n}(z,z')$ with different indices. The derivations for summations over indices in the polynomials $\,{\rm L}_{m,n}(z,z')$ (e.g., in generating functions) can be temporarily shifted in such way to derivations for the monomials $z^mz'^{\,n}$ with final application of the operator $\exp\left(-\frac{\partial^2}{\partial z \partial z'}\right)$ to the intermediate result.
In Section 11 we demonstrated the differences between both methods in one of the most simple cases which is the derivation of the well-known generating function (\ref{genherm}) for Hermite polynomials. However, our main aim was the derivation of new bilinear generating functions. In the bilinear generating functions for Hermite polynomials (\ref{mehler}) as well as for Laguerre 2D polynomials (\ref{gflag2dprfact}) we found interesting factorizations which establish connections to more special (linear) generating functions for these polynomials with transformed variables. Due to the rudimentary character of the generalized Laguerre (-Sonin) polynomials ${\rm L}_n^{\nu}(u),\,(u=zz')$ within the set of Laguerre 2D polynomials ${\rm L}_{m,n}(z,z')$ many formulae for the usual Laguerre polynomials become simpler and rigged with more symmetries if expressed by the Laguerre 2D polynomials.

From the bilinear generating functions or generating functions for products of Laguerre 2D polynomials one can derive completeness relations for the corresponding polynomials in the way such as demonstrated. Furthermore, we derived a simple sum (in the sense of not double sum!) over products of special Laguerre 2D polynomials which can be taken as intermediate step for the derivation of the generating function but this formula possesses proper importance for other calculations and was already useful in an application in quantum optics of phase states. The number of generating functions and of relations for Laguerre 2D and Hermite 2D polynomials is relatively large and a main source for suggestion are known generating functions and relations for usual Laguerre and, in particular, for Hermite polynomials.

The three generating functions for products of Hermite 2D and Laguerre 2D polynomials
(\ref{gflag2dpr}), (\ref{gfh2dpr}) and (\ref{gflagh4}) can be considered as special cases
of generating functions for products of Hermite 2D polynomials ${\rm H}_{m,n}({\sf
U};x,y){\rm H}_{m,n}({\sf V};u,v)$ or Laguerre 2D polynomials ${\rm L}_{m,n}({\sf
U};z,z'){\rm L}_{m,n}({\sf V};w',w)$ with arbitrary 2D matrices ${\sf U}$ and ${\sf V}$
as parameter mentioned in the Introduction. It is clear from the calculated different special cases that such generating functions are very complicated. Some simplification can be obtained by special choice of the matrix ${\sf V}$ related to the matrix ${\sf U}$, for example ${\sf V}^{\sf T}=\left({\sf U}^{-1}\right)$ where ${\sf A}^{\sf T}$ denotes the transposed matrix to ${\sf A}$. We began such calculations in \cite{w6} but to finish this is a task of future depending also on the appearance of problems in applications, for example, in quantum optics and classical optics (e.g., propagation of Gaussian beams) which require these generating functions.

We hope that we could convince the reader of some advantages of the use of Laguerre 2D polynomials in comparison to usual generalized Laguerre (-Sonin) polynomials as their radial rudiments. Although the usual Laguerre (and also Hermite) polynomials are mostly present in readily programmed form in mathematical computer programs it is not difficult to programme in the same way the Laguerre 2D polynomials by their explicit formulae (\ref{lag2d1}) as finite simple sums.

A main region of application of the derived generating functions is quantum optics of two harmonic oscillator modes and of quasiprobabilities of oscillator states such as the Wigner function and in classical optics the theory of Gauss-Hermite and Gauss-Laguerre beams and we applied some of the here derived relations in papers of former time 
[37--39].

In the following Appendix A, we develop some basic operator identities which are useful for calculation of convolutions of one- and two-dimensional Gaussian functions in combination
with polynomials and which were used in most of our derivations of the generating functions. The corresponding operators are connected with the Lie group $SU(1,1)$ (see, e.g., 
[40--43]).

\setcounter{chapter}{12}
\setcounter{equation}{0}
\aeqn
\section*{Appendix A:\\  Operator identities related to one- and two-dimensional Gaussian
convolutions by means of ${\cf {SU(1,1)}}$ operator disentanglement}

We sketch in this Appendix the derivation of some mostly novel and useful operational formulae
related to convolutions of Gaussian functions of one and two variables and use for this
purpose the technique of operator disentanglement of $SU(1,1)$ operators.

As canonical basis of the abstract Lie algebra to $SU(1,1)\sim SL(2,{\mathbb {R}})\sim Sp(2,{\mathbb {R}})$ and of its complex extension $SL(2,{\mathbb {C}})\sim Sp(2,{\mathbb {C}})$ are usually taken three
operators $(K_-,K_0,K_+)$ which obey the following commutation relations (e.g.,
[40--43]
and \cite{w7})
\begin{eqnarray}
\left[K_-,K_+\right]= 2K_0,\quad \left[K_0,K_-\right]=-K_-,\quad
\left[K_0,K_+\right]=+K_+.
\end{eqnarray}
As first case, we consider the following realization of the operators $(K_-,K_0,K_+)$ by
one-dimensional differentiation and multiplication operators
\begin{eqnarray}
K_-\equiv \frac{1}{2}\frac{\partial^2}{\partial x^2},\quad K_0\equiv
\frac{1}{4}\left(x\frac{\partial}{\partial x}+ \frac{\partial}{\partial x}x\right),\quad
K_+\equiv \frac{1}{2}x^2.\label{sureal1}
\end{eqnarray}
Using the commutation ($r$ is scalar parameter)
\begin{eqnarray}
\exp\left(\frac{x^2}{r}\right)\frac{\partial}{\partial
x}\exp\left(-\frac{x^2}{r}\right) &=& \frac{\partial}{\partial
x}+\frac{1}{1!r}\left[x^2,\frac{\partial}{\partial x}\right]+\frac{1}{2!r^2}\left[x^2,\left[x^2,\frac{\partial}{\partial x}\right]\right]+\ldots
\nonumber\\ &=& \frac{\partial}{\partial
x}-\frac{2}{r}x,
\end{eqnarray}
we obtain the following relation ($s$ is a second scalar parameter)
\begin{eqnarray}
\exp\left(\frac{s}{4}\frac{\partial^2}{\partial
x^2}\right)\exp\left(-\frac{x^2}{r}\right) &=& \exp\left(-\frac{x^2}{r}\right)\exp\left\{
\frac{s}{4}\left(\frac{\partial}{\partial x}-\frac{2}{r}x\right)^2\right\}\nonumber\\ &=&
\exp\left(-\frac{2}{r}K_+\right)\exp\left(\frac{s}{2}K_- -\frac{s}{r}2K_0
+\frac{2s}{r^2}K_+\right).\label{dis1}
\end{eqnarray}
Now, we can apply the following disentanglement relation for general group elements of
$Sp(2,{\mathbb {C}})\sim SL(2,{\mathbb {C}})$ (complexification of $SU(1,1)\sim Sp(2,{\mathbb {R}})\sim
SL(2,{\mathbb {R}})$) which is the first of the 6 relations with different ordering of the
factors derived in \cite{w7,w8}
\begin{eqnarray}
\exp\left(\xi K_- +\I \eta\, 2K_0 -\zeta K_+\right) =
\exp\left(-\frac{\mu}{\kappa}K_+\right)\exp\left(\lambda\kappa K_-\right)
\exp\Big(-\left(\log \kappa\right)2K_0\Big),\label{sudisent}
\end{eqnarray}
where $(\kappa,\lambda,\mu,\nu)$ are the matrix elements of the two-dimensional
fundamental representation of $Sp(2,{\mathbb {C}})$ in the basis of operators $(A,A^\dagger)$
(boson annihilation and creation operators $(a,a^\dagger)$ in simplest quantum-optical realization or $\left(\frac{\partial}{\partial x},x\right)$ in present case) forming together with
operators $(K_-,K_0,K_+)$ a basis of the Lie algebra to the inhomogeneous symplectic Lie
group $ISp(2,{\mathbb {C}})$. These elements which form an unimodular matrix (determinant equal to $1$) are explicitly
\begin{eqnarray}
\left(\begin{array}{cc}\kappa,& \lambda \\  \mu,& \nu
\end{array}\right) = \left(\begin{array}{cc}\dis{{\rm ch}(\varepsilon)
-\I \eta \frac{{\rm sh}(\varepsilon)}{\varepsilon}},&
\dis{\xi\frac{{\rm sh}(\varepsilon)}{\varepsilon}}\\
\dis{\zeta\frac{{\rm sh}(\varepsilon)}{\varepsilon}},& \dis{{\rm ch}(\varepsilon) + \I \eta
\frac{{\rm sh}(\varepsilon)}{\varepsilon}}
\end{array}\right),\quad \varepsilon \equiv \sqrt{\xi\zeta -\eta^2}.\label{sudisent1}
\end{eqnarray}
With the specialization $\xi=\frac{s}{2},\eta=\I\frac{s}{r},\zeta=-\frac{2s}{r^2}$, we
find the following specialization of this unimodular matrix (\ref{sudisent1})
\begin{eqnarray}
\left(\begin{array}{cc}\kappa,& \lambda \\  \mu,& \nu
\end{array}\right) = \left(\begin{array}{cc}\dis{1+\frac{s}{r},}&
\dis{\frac{s}{2}}\\[2mm]
\dis{-\frac{2s}{r^2},}& \dis{1-\frac{s}{r}}\end{array}\right),\quad \vep = 0,
\end{eqnarray}
and the disentanglement searched for is
\begin{eqnarray}
\exp\left(\frac{s}{2}K_- -\frac{s}{r}2K_0 +\frac{2s}{r^2}K_+\right) &=&
\exp\left(\frac{2s}{r(r+s)}K_+\right)\exp\left(\frac{s(r+s)}{2r}K_-\right)
\left(\frac{r}{r+s}\right)^{2K_0}.\qquad
\end{eqnarray}
Inserting this into (\ref{dis1}) and going back to realization (\ref{sureal1}), we obtain
the following important operator identity (we use $2K_0 \equiv\frac{1}{2}\left(x\frac{\partial}{\partial x}+ \frac{\partial}{\partial x}x\right) =x\frac{\partial}{\partial x} +\frac{1}{2})$
\begin{eqnarray}
\exp\left(\frac{s}{4}\frac{\partial^2}{\partial
x^2}\right)\exp\left(-\frac{x^2}{r}\right) = \sqrt{\frac{r}{r+s}}
\exp\left(-\frac{x^2}{r+s}\right)\exp\left(\frac{s(r+s)}{4r}\frac{\partial^2}{\partial
x^2}\right)\left(\frac{r}{r+s}\right)^{x\frac{\partial}{\partial x}}.\label{op1}
\end{eqnarray}
If we apply this operator identity to an arbitrary function $f(x)$, we can give this
relation a form which is often appropriate for direct application. Using that $x^n$ is an
eigenfunction of the operator $x\frac{\partial}{\partial x}$ to eigenvalue $n$, we find
\begin{eqnarray}
\left(x\frac{\partial}{\partial x}\right)x^n = nx^n, \quad\Rightarrow \quad\exp\left(\gamma x\frac{\partial}{\partial x}\right)x^n
=\left(\E^{\gamma}x\right)^n,\label{eigvaleq}
\end{eqnarray}
from which follows that $\exp\left(\gamma x\frac{\partial}{\partial x}\right)$ is the operator of multiplication of the argument of a function $f(x)$ according to \cite{w7}\footnote{The result of such and similar derivations (possibly with restriction to real $\gamma$) is also true for (possibly non-analytic) generalized functions $f(x)$ such as the step function $\theta(x)$ and the delta function and its derivatives $\delta^{(n)}(x)$ and many others since by their definition as linear continuous functionals, for example, the Taylor series expansion can be transformed to a sufficiently well-behaved class of basis functions and finally one can go back to the generalized functions.}
\begin{eqnarray}
\exp\left(\gamma x\frac{\partial}{\partial x}\right)f(x)=\exp\left(\gamma
x\frac{\partial}{\partial x}\right)\sum_{n=0}^\infty \frac{f^{(n)}(0)}{n!}x^n =
\sum_{n=0}^\infty \frac{f^{(n)}(0)}{n!}\left(\E^{\gamma}x\right)^n =
f\left(\E^{\gamma}x\right).\label{scalop}
\end{eqnarray}
By applying  the operator identity (\ref{op1}) to an arbitrary function $f(x)$ and using (\ref{scalop}) we obtain
\begin{eqnarray}
\exp\left(\frac{s}{4}\frac{\partial^2}{\partial
x^2}\right)\exp\left(-\frac{x^2}{r}\right)f(x) = \sqrt{\frac{r}{r+s}}
\exp\left(-\frac{x^2}{r+s}\right)\exp\left(\frac{s(r+s)}{4r}\frac{\partial^2}{\partial
x^2}\right)f\left(\frac{rx}{r+s}\right).\label{op2}
\end{eqnarray}
With substitution of variable $x \rightarrow y=\frac{r}{r+s}\,x$ and then by introduction of the new parameters $r'=\frac{r^2}{r+s},\,s'=\frac{rs}{r+s}$ or $r=r'+s',\,s=\frac{s'(r'+s')}{r'}$ we find from this relation
\begin{eqnarray}
\exp\left(-\frac{y^2}{r'}\right)\exp\left(\frac{s'}{4}\frac{\partial^2}{\partial
y^2}\right)f(y) = \sqrt{\frac{r'+s'}{r'}} \exp\left(\frac{r's'}{4(r'+s')}\frac{\partial^2}{\partial y^2}\right)\exp\left(-\frac{r'+s'}{r'^{2}}y^2\right)f\left(\frac{r'+s'}{r'}y\right).\nonumber\\
\end{eqnarray}
which is a sometimes useful transformation of (\ref{op2}).

Choosing the function $f(x)=\frac{1}{\sqrt{\pi r}}$, we obtain from (\ref{op2}) as special case the following formula
\begin{eqnarray}
\exp\left(\frac{s}{4}\frac{\partial^2}{\partial x^2}\right)\exp\left(-\frac{x^2}{r}\right)\frac{1}{\sqrt{\pi
r}} &\equiv & \frac{1}{\sqrt{\pi s}}\,\exp\left(-\frac{x^2}{s}\right)* \frac{1}{\sqrt{\pi r}}\,\exp\left(-\frac{x^2}{r}\right) \nonumber\\ &=& \frac{1}{\sqrt{\pi(r+s)}} \exp\left(-\frac{x^2}{r+s}\right),\label{convtwogauss}
\end{eqnarray}
which as shown may also be written as the convolution of two normalized Gaussian functions and provides as result again a normalized Gaussian function as it is well known (notation '$*$' means forming the convolution of two functions). The operator $\exp\left(\frac{s}{4}\frac{\partial^2}{\partial x^2}\right)$ applied to an arbitrary function makes the convolution of the normalized Gaussian function $\frac{1}{\sqrt{\pi s}}\exp\left(-\frac{x^2}{s}\right)$ with this arbitrary function that
can be proved by Fourier transformation. Below we find another derivation of this equivalence (see
(\ref{convgaussgf})). However, the main power of relation (\ref{op2}) is seen by
applying it to more complicated functions. For example, by applying it to the functions
$f(x)=(2x)^n$ and using the alternative definition of Hermite polynomials in first line
of (\ref{hermite}), we obtain
\begin{eqnarray}
&& \exp\left(\frac{s}{4}\frac{\partial^2}{\partial
x^2}\right)\exp\left(-\frac{x^2}{r}\right)(2x)^n \nonumber\\ &=& \sqrt{\frac{r}{r+s}}
\exp\left(-\frac{x^2}{r+s}\right)\left(\frac{\sqrt{-rs(r+s)}}{r+s}\,\right)^n {\rm
H}_n\left(\frac{rx}{\sqrt{-rs(r+s)}}\right),\label{op1hermpol}
\end{eqnarray}
where due to ${\rm H}_n(-x)=(-1)^n\,{\rm H}_n(x)$ one has to choose the same but arbitrary sign of the two roots $\sqrt{-rs(r+s)}$. The expression on the right-hand side of (\ref{op1hermpol}) is invariant with respect to the choice of this sign but if we shorten the fractions then the correlation of the signs of the now two different roots becomes unclear\footnote{This is the reason why we do not set $\frac{\sqrt{-rs(r+s)}}{r+s}=\sqrt{-\frac{rs}{r+s}}$.}.
In application to the Hermite polynomials ${\rm H}_n(x)$ we find in similar way
\begin{eqnarray}
&& \exp\left(\frac{s}{4}\frac{\partial^2}{\partial
x^2}\right)\exp\left(-\frac{x^2}{r}\right){\rm H}_n(x)\nonumber\\ &=&
\sqrt{\frac{r}{r+s}}
\exp\left(-\frac{x^2}{r+s}\right)\left(\frac{\sqrt{(r+s-rs)(r+s)}}{r+s}\,\right)^n {\rm
H}_n\left(\frac{rx}{\sqrt{(r+s-rs)(r+s)}}\,\right),\qquad \label{op2hermpol}
\end{eqnarray}
where again one can choose an arbitrary sign of the roots $\sqrt{(r+s-rs)(r+s)}$ which has only to be the same throughout the whole expression on the right-hand side that hinders us to shorten the formula.

As second case, we now consider the following realization of the operators
$(K_-,K_0,K_+)$ by two-dimensional differentiation and multiplication operators
\begin{eqnarray}
K_-\equiv \frac{\partial^2}{\partial z\partial z'},\quad K_0\equiv
\frac{1}{2}\left(z\frac{\partial}{\partial z}+\frac{\partial}{\partial z'}z'\right),\quad
K_+\equiv zz'.\label{sureal2}
\end{eqnarray}
Our notation of the two independent (in general, complex) variables $(z,z')$ is due to
the fact that in most potential applications of the formulae which we derive, we have a
pair of complex conjugated variables $(z,z^*)$ and we can then easily set $z'=z^*$ but
the results can also be applied if we have instead of this a pair of real variables
$(x,y)$. From the commutation
\begin{eqnarray}
\exp\left(\frac{zz'}{r}\right) \frac{\partial^2}{\partial z\partial
z'}\exp\left(-\frac{zz'}{r}\right)=\left(\frac{\partial}{\partial z}-\frac{z'}{r}\right)
\left(\frac{\partial}{\partial z'}-\frac{z}{r}\right),
\end{eqnarray}
follows
\begin{eqnarray}
\exp\left(s\frac{\partial^2}{\partial z\partial z'}\right)\exp\left(-\frac{zz'}{r}\right)
&=& \exp\left(-\frac{zz'}{r}\right)\exp\left\{s\left(\frac{\partial}{\partial
z}-\frac{z'}{r}\right) \left(\frac{\partial}{\partial z'}-\frac{z}{r}\right)\right\}
\nonumber\\ &=& \exp\left(-\frac{1}{r}K_+\right)\exp\left(sK_- -\frac{s}{r} 2K_0
+\frac{s}{r^2}K_+\right).\label{dis2}
\end{eqnarray}
The operator which we have to disentangle corresponds to the special choice
$\xi=s,\eta=\I\frac{s}{r},\zeta=-\frac{s}{r^2}$ in (\ref{sudisent}) and the unimodular
matrix (\ref{sudisent1}) takes on the following special form
\begin{eqnarray}
\left(\begin{array}{cc}\kappa,& \lambda \\  \mu,& \nu
\end{array}\right) = \left(\begin{array}{cc}\dis{1+\frac{s}{r},}&
\dis{s}\\[2mm]
\dis{-\frac{s}{r^2},}& \dis{1-\frac{s}{r}}\end{array}\right),\quad \vep \equiv \sqrt{\xi\zeta -\eta^2} = 0.
\end{eqnarray}
Using the disentanglement relation (\ref{sudisent}), from (\ref{dis2}) follows the
important operational identity
\begin{eqnarray}
&&\exp\left(s\frac{\partial^2}{\partial z\partial
z'}\right)\exp\left(-\frac{zz'}{r}\right) \nonumber\\ &=&
\frac{r}{r+s}\exp\left(-\frac{zz'}{r+s}\right) \exp\left(\frac{s(r+s)}{r}
\frac{\partial^2}{\partial z\partial z'}\right) \left(\frac{r}{r+s}
\right)^{z\frac{\partial}{\partial z}+z'\frac{\partial}{\partial z'}},\label{op3}
\end{eqnarray}
Applied to arbitrary functions $f(z,z')$, we find similar to (\ref{op2})
\begin{eqnarray}
&&\exp\left(s\frac{\partial^2}{\partial z\partial
z'}\right)\exp\left(-\frac{zz'}{r}\right)f(z,z')\nonumber\\ &=&
\frac{r}{r+s}\exp\left(-\frac{zz'}{r+s}\right) \exp\left(\frac{s(r+s)}{r}
\frac{\partial^2}{\partial z\partial z'}\right)f\left(\frac{r z}{r+s},\frac{r
z'}{r+s}\right),\label{op4}
\end{eqnarray}
which in special case $f(z,z')=1$ possesses some relation to the generating function for usual Laguerre polynomials (see (\ref{convtwogauss}) below and (\ref{genlagpol}) and (\ref{genlagpol1})).
By transition to new variables $z\rightarrow w=\frac{r}{r+s}\,z,\;z'\rightarrow w'=\frac{r }{r+s}\,z'$ in (\ref{op4}) and changing the parameters to $r'=\frac{r^2}{r+s}\,,s'=\frac{rs}{r+s}$ or $r=r'+s',\,s=\frac{s'(r'+s')}{r'}$ we obtain
\begin{eqnarray}
&& \exp\left(-\frac{ww'}{r'}\right)\exp\left(s'\frac{\partial^2}{\partial w\partial
w'}\right)f(w,w') \nonumber\\ &=& \frac{r'+s'}{r'}\exp\left(\frac{r's'}{r'+s'}\frac{\partial^2}{\partial w\partial w'}\right)\exp\left(-\frac{r'+s'}{r'^2}ww'\right)f\left(\frac{r'+s'}{r'}w,\frac{r'+s'}{r'}w'\right),
\end{eqnarray}
as some useful transformation of (\ref{op4}).

We mention that the monomials $z^mz'^n$ are eigenfunctions of the operator $2K_0$ in (\ref{sureal2}) to the eigenvalues $ m+n+1$ according to
\begin{eqnarray}
\left(z\frac{\partial}{\partial z}+\frac{\partial}{\partial z'}z'\right)z^mz'^n &=& (m+n+1)z^mz'^n,
\end{eqnarray}
from which follows in application of its exponential to an arbitrary function $f(z,z')$
\begin{eqnarray}
\exp\left(\lambda\left(z\frac{\partial}{\partial z}+\frac{\partial}{\partial z'}z'\right)\right)f(z,z')
&=& \E^{\lambda}f\left(\E^{\lambda} z,\E^{\lambda}z'\right),
\end{eqnarray}
in analogy to formula (\ref{scalop}) together with (\ref{eigvaleq}). This was used in (\ref{op4}).

In special case of function $f(z,z')= \frac{1}{\pi r}$, we obtain from (\ref{op4}) the two-dimensional convolution of two (normalized if $z'=z^*)$ Gaussian functions with parameters $r$ and $s$
\begin{eqnarray}
\exp\left(s\frac{\partial^2}{\partial z\partial
z'}\right)\exp\left(-\frac{zz'}{r}\right)\frac{1}{\pi r} &=& \frac{1}{\pi s}\,\exp\left(-\frac{zz'}{s}\right)* \frac{1}{\pi r}\,\exp\left(-\frac{zz'}{r}\right)\nonumber\\
&=& \frac{1}{\pi(r+s)}\exp\left(-\frac{zz'}{r+s}\right),\label{convtwogauss}
\end{eqnarray}
which provides again a normalized Gaussian function with the parameter $r+s$.

The full power of (\ref{op4}) is seen if we apply it to the functions $f(z,z')=z^mz'^n$
and if we use the definition of Laguerre 2D polynomials given in first line of
(\ref{lag2d1}) that leads to
\begin{eqnarray}
&& \exp\left(s\frac{\partial^2}{\partial z\partial
z'}\right)\exp\left(-\frac{zz'}{r}\right)z^{m}z'^{\,n} \;=\; \frac{r}{r+s}
\exp\left(-\frac{zz'}{r+s}\right)\nonumber\\
&& \hspace{15mm} \cdot \left(\frac{\sqrt{-rs(r+s)}}{r+s}\right)^{m+n}{\rm
L}_{m,n}\left(\frac{rz}{\sqrt{-rs(r+s)}},\frac{rz'}{\sqrt{-rs(r+s)}}\right),
\end{eqnarray}
Applied to Laguerre 2D polynomials $f(z,z')={\rm L}_{m,n}(z,z')$, we find in analogous
way
\begin{eqnarray}
&& \exp\left(s\frac{\partial^2}{\partial z\partial
z'}\right)\exp\left(-\frac{zz'}{r}\right){\rm L}_{m,n}(z,z') \;=\; \frac{r}{r+s}
\exp\left(-\frac{zz'}{r+s}\right) \nonumber\\
&& \hspace{2mm} \cdot\left(\frac{\sqrt{(r+s-rs)(r+s)}}{r+s}\,\right)^{m+n}{\rm
L}_{m,n}\left(\frac{rz}{\sqrt{(r+s-rs)(r+s)}},
\frac{rz'}{\sqrt{(r+s-rs)(r+s)}}\right).\qquad
\end{eqnarray}
In last two relations one has to choose an arbitrary but the same sign within one formula for the roots $\sqrt{-rs(r+s)}$ or $\sqrt{(r+s-rs)(r+s)}$, respectively, due to ${\rm L}_{m,n}(-z,-z')= (-1)^{m+n}\,{\rm L}_{m,n}(z,z')$ that means for reason that was explained already for the analogous formulae (\ref{op1hermpol}) and (\ref{op2hermpol}).

We mention that an additional displacement of the arguments does not make any
difficulties in the derived operator relations. For example, according to
\begin{eqnarray}
\exp\left(-x_0\frac{\partial}{\partial x}\right)x\exp\left(x_0\frac{\partial}{\partial
x}\right)=x-x_0,\;\;\Rightarrow \;\; f(x-x_0)=\exp\left(-x_0\frac{\partial}{\partial
x}\right)f(x)\exp\left(x_0\frac{\partial}{\partial x}\right),\qquad \label{displop}
\end{eqnarray}
we can generalize the operator identity (\ref{op1}) in the following way
\begin{eqnarray}
&&\exp\left(\frac{s}{4}\frac{\partial^2}{\partial
x^2}\right)\exp\left(-\frac{\left(x-x_0\right)^2}{r}\right) \nonumber\\ &=&
\sqrt{\frac{r}{r+s}}
\exp\left(-\frac{\left(x-x_0\right)^2}{r+s}\right)\exp\left(\frac{s(r+s)}{4r}
\frac{\partial^2}{\partial x^2}\right)\left(\frac{r}{r+s}\right)^{\left(x-x_0\right)
\frac{\partial}{\partial x}},\label{displop1}
\end{eqnarray}
and in application to an arbitrary function $f(x)=f\left(x_0+x-x_0\right)$
\begin{eqnarray}
&&\exp\left(\frac{s}{4}\frac{\partial^2}{\partial
x^2}\right)\exp\left(-\frac{\left(x-x_0\right)^2}{r}\right)f(x) \nonumber\\ &=&
\sqrt{\frac{r}{r+s}}
\exp\left(-\frac{\left(x-x_0\right)^2}{r+s}\right)\exp\left(\frac{s(r+s)}{4r}
\frac{\partial^2}{\partial x^2}\right)f\left(\frac{rx+sx_0}{r+s}\right).\label{displop2}
\end{eqnarray}
Choosing $f(x)=\frac{1}{\sqrt{\pi r}}$ in (\ref{displop2}) and making then the limiting transition $r\rightarrow 0$ using
\begin{eqnarray}
\lim_{r\rightarrow 0}\frac{1}{\sqrt{\pi r}}\exp\left(-\frac{\left(x-x_0\right)^2}{r}\right) = \delta\left(x-x_0\right), \label{onedeltafunc}
\end{eqnarray}
one finds the transformation
\begin{eqnarray}
\exp\left(\frac{s}{4}\frac{\partial^2}{\partial x^2}\right)\delta\left(x-x_0\right) &=& \frac{1}{\sqrt{\pi s}}\exp\left(-\frac{\left(x-x_0\right)^2}{s}\right),
\end{eqnarray}
and by multiplication of both sides with an arbitrary function $g(x_0)$ and by integration over $x_0$
\begin{eqnarray}
\exp\left(\frac{s}{4}\frac{\partial^2}{\partial x^2}\right)g(x) &=& \frac{1}{\sqrt{\pi s}} \int_{-\infty}^{+\infty}dx_0\,\exp\left(-\frac{\left(x-x_0\right)^2}{s}\right)g(x_0) \nonumber\\ &\equiv & \frac{1}{\sqrt{\pi s}}\exp\left(-\frac{x^2}{s}\right)*g(x),\label{convgaussgf}
\end{eqnarray}
that are two different representations of the convolution of a normalized Gaussian function with an arbitrary function $g(x)$. A special case is the convolution of two normalized Gaussian functions with the result given in (\ref{convtwogauss}). For $s\le 0$ the convolutions have to be considered in the sense of the theory of generalized functions.

Analogously to (\ref{displop}), relation (\ref{op3}) can be generalized including two, in general, different complex displacements $z\rightarrow z-z_0$ and $z'\rightarrow z'-z'_0$ leading to the operator identity
\begin{eqnarray}
&&\exp\left(s\frac{\partial^2}{\partial z\partial
z'}\right)\exp\left(-\frac{\left(z-z_0\right)\left(z'-z'_0\right)}{r}\right) \nonumber\\ &=&
\frac{r}{r+s}\exp\left(-\frac{\left(z-z_0\right)\left(z'-z'_0\right)}{r+s}\right) \exp\left(\frac{s(r+s)}{r} \frac{\partial^2}{\partial z\partial z'}\right) \left(\frac{r}{r+s}
\right)^{\left(z-z_0\right)\frac{\partial}{\partial z}+\left(z'-z'_0\right)\frac{\partial}{\partial z'}},\qquad\label{displop3}
\end{eqnarray}
with the consequence
\begin{eqnarray}
&&\exp\left(s\frac{\partial^2}{\partial z\partial
z'}\right)\exp\left(-\frac{\left(z-z_0\right)\left(z'-z'_0\right)}{r}\right)f(z,z')\nonumber\\ &=&
\frac{r}{r+s}\exp\left(-\frac{\left(z-z_0\right)\left(z'-z'_0\right)}{r+s}\right) \exp\left(\frac{s(r+s)}{r}\frac{\partial^2}{\partial z\partial z'}\right)f\left(\frac{r z+sz_0}{r+s},\frac{rz'+sz'_0}{r+s}\right).\qquad\label{displop4}
\end{eqnarray}
By a limiting procedure in analogy to the derivation of (\ref{convgaussgf}) one obtains from (\ref{displop4}) for $z'=z^*$ where $z^*$ is complex conjugate to $z$
\begin{eqnarray}
\exp\left(s\frac{\partial^2}{\partial z\partial z^*}\right)g\left(z,z^*\right) &=& \frac{1}{\pi s}\int \frac{\I}{2}dz_0\wedge dz^*_0 \exp\left(-\frac{\left(z-z_0\right)\left(z^*-z^*_0\right)}{s}\right)g\left(z_0,z^*_0\right) \nonumber\\ &\equiv & \frac{1}{\pi s}\exp\left(-\frac{zz^*}{s}\right)*g\left(z,z^*\right),
\end{eqnarray}
which are two different representations of the two-dimensional convolution of a normalized Gaussian function with an arbitrary function, in particular
\begin{eqnarray}
\exp\left(s\frac{\partial^2}{\partial z\partial z^*}\right)\delta\left(z-z_0,z^*-z_0^*\right) &=& \frac{1}{\pi s}\exp\left(-\frac{\left(z-z_0\right)\left(z^*-z_0^*\right)}{s}\right),
\end{eqnarray}
with the two-dimensional delta function obtained by the limiting procedure
\begin{eqnarray}
\lim_{r\rightarrow 0}\frac{1}{\pi r}\exp\left(-\frac{\left(z-z_0\right)\left(z^*-z_0^*\right)}{r}\right) = \delta\left(z-z_0,z^*-z_0^*\right), \label{twodeltafunc}
\end{eqnarray}
in analogy to (\ref{onedeltafunc}).
All this can be also obtained by Fourier transformation and its inversion.

The displacement operation (\ref{displop}) with substitutions $x \rightarrow z,\frac{\partial}{\partial x}\rightarrow \frac{\partial}{\partial z}$ and with the formal substitution $x_0 \rightarrow \frac{\partial}{\partial z'}$ can also be applied to evaluate the following expression
\begin{eqnarray}
\exp\left(-\frac{\partial^2}{\partial z\partial
z'}\right)\exp\left(-\frac{\sigma^2}{2}z^2-\frac{\tau^2}{2}z'^2\right) &=&
\exp\left\{-\frac{\sigma^2}{2}\left(z-\frac{\partial}{\partial z'}\right)^2\right\}
\exp\left(-\frac{\tau^2}{2}z'^2\right)\nonumber\\ &=&
\frac{1}{\sqrt{1-\sigma^2\tau^2}}\exp\left(-\frac{\sigma^2z^2+\tau^2z'^2+2\sigma^2\tau^2 zz'}
{2\left(1-\sigma^2\tau^2\right)}\right).\qquad \label{expdispl}
\end{eqnarray}
Necessarily, this agrees with the specialized generating function (\ref{gflagh5}). We
used this formula for the evaluation of (\ref{gflagh2}) in one of some possible variants of proof.
In application of the operator $\exp\left(-\frac{\partial^2}{\partial z\partial
z'}\right)$ to an arbitrary function $f(z,z')$ we find as generalization of (\ref{expdispl})
\begin{eqnarray}
\exp\left(s\frac{\partial^2}{\partial z\partial
z'}\right) f(z,z') &=& f\left(z+s\frac{\partial}{\partial z'},z'+s\frac{\partial}{\partial z}\right)1 \nonumber\\ &=& \sum_{m=0}^\infty\sum_{n=0}^\infty \frac{f^{(m,n)}(0,0)}{m!n!}\left(\sqrt{-s}\,\right)^{m+n}\,{\rm L}_{m,n}\left(\frac{z}{\sqrt{-s}},\frac{z'}{\sqrt{-s}}\right),
\end{eqnarray}
where we applied the Taylor series of $f(z,z')$ in powers of $z$ and $z'$ and the definition (\ref{lag2d1}) of the Laguerre 2D polynomials and where we took into account that the operators $z+s\frac{\partial}{\partial z'}$ and $z'+s\frac{\partial}{\partial z}$ are commuting. For example, if we specialize the function to $f(z,z')=\exp\left(-\frac{\sigma^2}{2}z^2-\frac{\tau^2}{2}z'^2\right)$ then we find according to (\ref{expdispl}) the identity
\begin{eqnarray}
&&\exp\left(s\frac{\partial^2}{\partial z\partial
z'}\right)\exp\left(-\frac{\sigma^2}{2}z^2\right)\exp\left(-\frac{\tau^2}{2}z'^2\right) \nonumber\\ &=&\sum_{k=0}^\infty\sum_{l=0}^\infty \frac{1}{k!l!}\left(\frac{s}{2}\right)^{k+l}\sigma^{2k}\tau^{2l}\,{\rm L}_{2k,2l}\left(\frac{z}{\sqrt{-s}},\frac{z'}{\sqrt{-s}}\right) \nonumber\\ &=& \frac{1}{\sqrt{1-s^2\sigma^2\tau^2}}\exp\left(-\frac{\sigma^2z^2+\tau^2z'^2-2s\sigma^2\tau^2 zz'}{2\left(1-s^2\sigma^2\tau^2\right)}\right),\qquad
\end{eqnarray}
which for $s=-1$ is identical with (\ref{gflagh5}).

The derived identities are particularly important in connection with the operational definition of the Hermite polynomials in (\ref{hermite}) and of the Laguerre 2D polynomials in (\ref{lag2d1}). It seems to us that the here chosen derivations belong to the most simple and perspective ones.

{\small{

}}

\end{document}